\documentclass[final,onefignum,onetabnum]{siamart220329}

\usepackage{amsfonts}
\usepackage{graphicx}
\usepackage{subcaption}
\usepackage{epstopdf}
\usepackage{algorithmic}
\usepackage{eqparbox}

\ifpdf
  \DeclareGraphicsExtensions{.eps,.pdf,.png,.jpg}
\else
  \DeclareGraphicsExtensions{.eps}
\fi

\newsiamremark{remark}{Remark}
\newsiamremark{hypothesis}{Hypothesis}
\crefname{hypothesis}{Hypothesis}{Hypotheses}
\newsiamthm{claim}{Claim}

\headers{An Efficient Scaled spectral preconditioner for linear systems}{Y. Diouane, S. G\"urol, O. Mouhtal and D. Orban}

\title{An Efficient Scaled spectral preconditioner for sequences of symmetric positive definite linear systems\thanks{Submitted to the editors DATE.
\funding{This work was funded by French National Programme LEFE/INSU.}}}

\author{Y. DIOUANE\thanks{GERAD and Department of Mathematics and Industrial Engineering, Polytechnique Montr\'eal. 
(\email{youssef.diouane@polymtl.ca},
\email{dominique.orban@polymtl.ca}).}
\and S.G\"urol\thanks{CERFACS / CECI
CNRS UMR 5318, Toulouse, France.
  (\email{gurol@cerfacs.fr}, \email{mouhtal@cerfacs.fr}).}
\and  O.MOUHTAL\footnotemark[3] \footnotemark[2]  
\and D.ORBAN\footnotemark[2]  
}

\usepackage{amsopn}

\usepackage[normalem]{ulem}
\usepackage{mathtools}
\usepackage{dsfont}
\usepackage{comment}

\newcommand{\R}{\mathds{R}}

\newcommand{\thetares}{\theta_{r}}
\newcommand{\thetamidrange}{\theta_{m}}

\ifpdf
\hypersetup{
  pdftitle={An Efficient Scaled spectral preconditioner for sequences of symmetric positive definite linear systems},
  pdfauthor={Youssef Diouane, Selime G\"urol, Oussama Mouhtal and Dominique Orban}
}
\fi

\begin{document}

\maketitle

\pagestyle{myheadings}

\begin{abstract}
We explore a \textit{scaled} spectral preconditioner for the efficient solution of sequences of symmetric and positive-definite linear systems. We design the scaled preconditioner not only as an approximation of the inverse of the linear system but also with consideration of its use within the conjugate gradient (CG) method. We propose three different strategies for selecting a scaling parameter, which aims to position the eigenvalues of the preconditioned matrix in a way that reduces the energy norm of the error, the quantity that CG monotonically decreases at each iteration. Our focus is on accelerating convergence especially in the early iterations, which is particularly important when CG is truncated due to computational cost constraints. Numerical experiments provide in data assimilation confirm that the \textit{scaled} spectral preconditioner can significantly improve early CG convergence with negligible computational cost.
\end{abstract}

\begin{keywords}
Sequence of linear systems, conjugate gradient method, deflated CG, spectral preconditioner, convergence rate, data assimilation.
\end{keywords}

\begin{MSCcodes}
68Q25, 68R10, 68U05
\end{MSCcodes}

\section{Introduction}

Efficiently solving sequences of symmetric positive-definite (SPD) linear systems
\begin{equation}%
    \label{eq:spd-sequence}
    A^{(j)} x^{(j)} = b^{(j)}, \quad  j = 1, 2, \ldots
\end{equation}
is crucial in various inverse problems of computational science and engineering. For instance, in data assimilation~\cite{Dale91, Kalnay2002}, where one aims to solve a large-scale weighted regularized nonlinear least-squares problem via the truncated Gauss-Newton algorithm (GN)~\cite{gratton2007approximate, NoceWrig06}, each iteration involves solving a linear least-squares subproblem. The latter may be formulated as a large SPD linear system, typically solved using the preconditioned conjugate-gradient method (PCG). Since consecutive systems do not differ significantly, recycling Krylov subspace information has been explored and proven to be effective~\cite{DataAssimilaton, AutoPrec, LMP, nash1991numerical}. 

One way of recycling Krylov subspace information involves leveraging search directions obtained from PCG on earlier systems to construct a limited-memory quasi-Newton preconditioner (LMP)~\cite{AutoPrec, nash1991numerical}.  This preconditioner, built solely from PCG information, does not require explicit knowledge of any matrix in the sequence, making it particularly suitable for applications where only matrix-vector products are available, which is the case of data assimilation.~\cite{LMP} generalizes this limited-memory preconditioner, and introduces specific variants when used with eigen- or Ritz pairs.

They focused on a first-level preconditioner, capable of clustering most eigenvalues at $1$ with few outliers, is already available for the first linear system in sequence. Then, they used LMP as a second-level preconditioner to improve the efficiency of the first. The goal of the LMP is to capture directions in a low-dimensional subspace that the first-level preconditioner may miss, and use them to improve convergence of PCG. When \(A^{(j)} = A\) for all \(j\), spectral analysis of the preconditioned matrix when used with $k$ pairs has shown that LMP can cluster at least $k$ eigenvalues at $1$, and that the eigenvalues of the preconditioned matrix interlace with those of the original matrix~\cite{LMP}. The efficiency of this approach has been demonstrated in a real-life data assimilation applications~\cite{LMP, tshimanga2007}. 

We focus on improving the performance of the \emph{spectral LMP}~\cite{SenSpectLMP, LMP}, which is built by using eigenpairs of \(A^{(j)}\). The spectral LMP shares the same formulation as the abstract balancing domain decomposition method~\cite{nabben2006comparison} and is equivalent to deflation-based preconditioning when used with a specific initial point~\cite{tshimanga2007}. 

When designing preconditioners for PCG, the primary focus in the literature is mostly on $A$ and the significance of the initial guess is overlooked. Although the importance of the initial guess is mentioned, its impact on the choice of a preconditioner is not well studied. Favorable eigenvalue distributions are also highlighted in terms of clustering, but there is little emphasis on the position of the clusters. The performance of the preconditioner is also measured in terms of the total number of iterations to converge, with little focus on the convergence in the early iterations. When PCG is truncated before convergence due to computational budget or when used as a solver within a optimization method like GN, the effect of the preconditioner on the early convergence of PCG is also crucial. In this paper, we aim to explore those overlooked aspects to design a good preconditioner. We not only aim to improve convergence by reducing the total number of iterations but also ensure that, from the very first iteration, the preconditioned iterates outperform those of the original system. In doing so, we specifically focus on strategically positioning the eigenvalues captured by the LMP, in that the energy norm of the error at each iteration of CG is reduced.

The paper is organized as follows. In \Cref{sec:notation} we start by introducing the necessary notation. In \Cref{sec:Background}, we review PCG and its convergence properties. We then discuss the characteristics of an efficient preconditioner that can be applied to~\eqref{eq:spd-sequence}. \Cref{sec:A scaled spectral preconditioner} is our main contribution. We define the \emph{scaled} spectral preconditioner and discuss its properties. Next, we outline three key approaches for selecting the scaling parameter, which influences the positioning of the eigenvalue cluster, to reduce total number of iterations and enhance convergence in the early iterations. In \Cref{sec:NumExp}, we provide numerical experiments using the Lorenz~\(95\) reference model from data assimilation to validate theoretical results. Finally, conclusions and perspectives are discussed in \Cref{sec:conclusion}.

\section{Notation}
\label{sec:notation}
The matrix $A \in \R^{n \times n}$ is always SPD.
Its spectral radius is $\rho(A)$.
Its spectral decomposition is $A = S \Lambda S^\top$ with $\Lambda = \text{diag}(\lambda_1, \ldots, \lambda_n)$, $ \lambda_1 \geq \ldots \geq \lambda_n > 0$, and $S = \begin{bmatrix} s_1 & \cdots & s_n \end{bmatrix}$ orthogonal.
Its $i$-th eigenvalue is  $\nu_i(A)$.
Its range space is $\mathcal{R}(A)$.
The $A$-norm, or \emph{energy norm}, of vector $x$ is $\|x\|_{A} =\sqrt{x^\top A x}$.
The spectral norm is $\| . \|_2$.

\section{Background}
\label{sec:Background}
\subsection{CG algorithm}
The Conjugate Gradient (CG) method~\cite{Hestenes1952MethodsOC} is the workhorse for $Ax = b$ with SPD $A \in \R^{n \times n}$ and $b \in \R^{n}$. If $x_0 \in \R^{n}$ is an initial guess and $r_0 = b - A x_0$ is the initial residual, then at every step $\ell = 1, 2, \ldots, n$, CG produces a unique approximation~\cite[p.176]{YoussSaad}
\begin{equation}
    x_\ell \in x_0 + \mathcal{K}_\ell(A, r_0) \quad \text{such that} \quad r_\ell \perp \mathcal{K}_\ell(A, r_0),\label{eq:projProcess}
\end{equation}
which is equivalent~\cite[p.126]{YoussSaad} to
\begin{equation}
     \|x^\ast - x_\ell\|_{A} = \min_{x \in x_0 + \mathcal{K}_\ell(A, r_0)} \left\| x^\ast - x \right\|_{A},\label{eq:MiniPro}
\end{equation}
where $x^\ast$ is the exact solution, $\mathcal{K}_\ell(A, r_0):= \text{span}\{r_0, A r_0, \ldots, A^{\ell-1} r_0\}$ is the $\ell$-th Krylov subspace generated by $A$ and $r_0$. 
In exact arithmetic, the method terminates in at most \(\mu\) iterations, where $\mu$ is the grade of $r_0$ with respect to $A$, i.e., the maximum dimension of the Krylov subspace generated by $A$ and $r_0$~\cite{YoussSaad}.
The most popular and computationally efficient variant of~\eqref{eq:projProcess} is the original formulation of~\cite{Hestenes1952MethodsOC}, that recursively updates coupled 2-term recurrences for $x_{\ell+1}$, $r_{\ell+1}$, and the search direction ${p}_{\ell+1}$. \Cref{Algo:CG} states the complete algorithm. A common stopping criterion is based on sufficient decrease of the relative residual norm. However, in practical data assimilation implementations, a fixed number of iterations is used as stopping criterion due to computational budget constraints.
CG is presented alongside its companion formulation, \Cref{Algo:PCG}, to be detailed in \Cref{subsec:Properties of a good preconditioner}.

\noindent
\begin{minipage}[t]{0.47\linewidth}
\begin{algorithm}[H]
\caption{CG\vphantom{PCG}}
\label{Algo:CG}
\begin{algorithmic}[1]
\STATE $r_0 = b - Ax_0$                               $\vphantom{\hat{r}_0 = b - A \hat{x}_0}$
\STATE $\vphantom{z_0 = F \hat{r}_0}$
\STATE $\rho_0 = r_0^\top r_0$                   $\vphantom{\hat{p}_0 = z_0}$
\STATE $p_0 = r_0$                               $\vphantom{\hat{p}_0 = z_0}$
\FOR{$\ell = 0, 1, \ldots$}
    \STATE $q_\ell = Ap_\ell$                               $\vphantom{\hat{q}_\ell = A \hat{p}_\ell}$
    \STATE $\alpha_\ell = \rho_\ell / (q_\ell^\top p_\ell)$       $\vphantom{\hat{\alpha}_\ell = \hat{\rho}_\ell / (\hat{q}_\ell^\top \hat{p}_\ell)}$
    \STATE $x_{\ell+1} = x_\ell + \alpha_\ell p_\ell$             $\vphantom{\hat{x}_{\ell+1} = \hat{x}_\ell + \hat{\alpha}_\ell \hat{p}_\ell}$
    \STATE $r_{\ell+1} = r_\ell - \alpha_\ell q_\ell$             $\vphantom{\hat{r}_{\ell+1} = \hat{r}_\ell - \hat{\alpha}_\ell \hat{q}_\ell}$                                           \STATE $\vphantom{z_{\ell+1} = F \hat{r}_{\ell+1}}$
    \STATE $\rho_{\ell+1} = r_{\ell+1}^\top r_{\ell+1}$   $\vphantom{\hat{\rho}_{\ell+1} = \hat{r}_{\ell+1}^\top z_{\ell+1}}$
    \STATE $\beta_{\ell+1} = \rho_{\ell+1} / \rho_\ell$        $\vphantom{\hat{\beta}_{\ell+1} = \hat{\rho}_{\ell+1} / \hat{\rho}_\ell}$    \STATE $p_{\ell+1} = r_{\ell+1} + \beta_{\ell+1} p_\ell$ $\vphantom{\hat{p}_{\ell+1} = z_{\ell+1} + \hat{\beta}_{\ell+1} \hat{p}_\ell}$
\ENDFOR
\end{algorithmic}
\end{algorithm}
\end{minipage}
\hfill
\begin{minipage}[t]{0.47\linewidth}
\begin{algorithm}[H]
\caption{PCG}
\label{Algo:PCG}
\begin{algorithmic}[1]
\STATE $\hat{r}_0 = b - A \hat{x}_0$ 
\STATE $z_0 = F \hat{r}_0$
\STATE $\hat{\rho}_0 = \hat{r}_0^\top z_0$
\STATE $\hat{p}_0 = z_0$
\FOR{$\ell = 0, 1, \ldots$}
    \STATE $\hat{q}_\ell = A \hat{p}_\ell$
    \STATE $\hat{\alpha}_\ell = \hat{\rho}_\ell / (\hat{q}_\ell^\top \hat{p}_\ell)$
    \STATE $\hat{x}_{\ell+1} = \hat{x}_\ell + \hat{\alpha}_\ell \hat{p}_\ell$ 
    \STATE $\hat{r}_{\ell+1} = \hat{r}_\ell - \hat{\alpha}_\ell \hat{q}_\ell$
    \STATE $z_{\ell+1} = F \hat{r}_{\ell+1}$
    \STATE $\hat{\rho}_{\ell+1} = \hat{r}_{\ell+1}^\top z_{\ell+1}$
    \STATE $\hat{\beta}_{\ell+1} = \hat{\rho}_{\ell+1} / \hat{\rho}_\ell$
    \STATE $\hat{p}_{\ell+1} = z_{\ell+1} + \hat{\beta}_{\ell+1} \hat{p}_\ell$
\ENDFOR
\end{algorithmic}
\end{algorithm}
\end{minipage}

\subsection{Convergence properties of CG}
 
The approximation $x_\ell$ uniquely determined by~\eqref{eq:projProcess} minimizes the error in the energy norm:
\begin{equation}
\|x^\ast - x_\ell\|_{A}^2 = \min_{p\in \mathbb{P}_\ell(0)} \|p(A)(x^\ast - x_0)\|_{A}^2 = \min_{p\in \mathbb{P}_\ell(0)} \, \sum_{i=1}^n p(\lambda_i)^2 \frac{\eta_i^2 }{\lambda_i}, \label{eq:MiniProcess}
\end{equation}
where $\eta_i = s_i^\top r_0$ and $\mathbb{P}_\ell(0)$ is the set of polynomials of degree at most $\ell$ with value $1$ at zero~\cite[p.193]{YoussSaad}. Thus, at each iteration, CG solves a certain weighted polynomial approximation problem over the discrete set~$\{\lambda_1, \ldots, \lambda_n\}$.
Moreover, if $z_1^{(\ell)}, \ldots, z_\ell^{(\ell)}$ are the $\ell$ roots of the solution $p_\ell^\ast$ to~\eqref{eq:MiniProcess},
\begin{equation}
\|x^\ast - x_\ell\|_{A}^2 = \sum_{i=1}^n p_\ell^\ast(\lambda_i)^2 \frac{\eta_i^2}{\lambda_i} = \sum_{i=1}^n \prod_{j=1}^\ell \left(1 - \frac{\lambda_i}{z_j^{(\ell)}}\right)^2 \frac{\eta_i^2}{\lambda_i}. \label{eq:ErrorRitz}
\end{equation}
The $z_j^{(\ell)}$ are the \emph{Ritz values}~\cite{VanDer}. From~\eqref{eq:ErrorRitz}, if $z_j^{(\ell)}$ is close to a $\lambda_i$, we expect a significant reduction in the error in energy norm.
Based on the above,~\cite{VanDer} explains the rate of convergence of CG in terms of the convergence of the Ritz values to eigenvalues of \(A\). Assuming that $\lambda_1, \ldots, \lambda_n$ take on the $r$ distinct values $\rho_1, \ldots, \rho_r$, 
CG converges in at most \(r\) iterations~\cite[Theorem 5.4]{NoceWrig06}. 

Using~\eqref{eq:MiniProcess} and maximizing over the values $p(\lambda_i)$~\cite[p.194]{YoussSaad}
leads to
\begin{equation}\label{eq:upper_bound2}
    \frac{\|x^\ast - x_\ell\|_A}{\|x^\ast - x_0\|_A} \leq \min_{p \in \mathbb{P}_\ell(0)} \max_{1 \leq i \leq n} |p(\lambda_i)|.
\end{equation}
By replacing $\{\lambda_1, \ldots , \lambda_n\}$ with the interval $[\lambda_1, \lambda_n]$ and using Chebyshev polynomials, we obtain an upper bound~\cite[p.194]{YoussSaad}:
\begin{equation}\label{eq:overestimate}
    \frac{\|x^\ast - x_\ell\|_A}{\|x^\ast - x_0\|_A} \leq 2 \left(\frac{\sqrt{\kappa(A)} - 1}{\sqrt{\kappa(A)} + 1} \right)^\ell,
\end{equation}
where $\kappa(A) := \lambda_1 / \lambda_n$ is the condition number of \(A\). While~\eqref{eq:upper_bound2} and~\eqref{eq:overestimate} provide the worst-case behavior of CG~\cite{PolyComp}, the convergence properties may vary significantly from the worst case for a specific initial approximation. Note also that upper bounds~\eqref{eq:upper_bound2} and~\eqref{eq:overestimate} only depend on A, and not on $r_0$. Though~\eqref{eq:overestimate} relates the convergence behavior of CG to \(\kappa(A)\), one should be careful as convergence is also influenced by the clustering of the eigenvalues and their positioning~\cite{carson2024towards, carson2020cost}.

\subsection{Properties of a good preconditioner}\label{subsec:Properties of a good preconditioner}
In many practical applications, a preconditioner is essential for accelerating the convergence of CG~\cite{Preconditionning1, wathen2015preconditioning}. Assume that a preconditioner $F = U U^\top \in \mathbb{R}^{n \times n}$ is available in a factored form, where $U$ is SPD, and consider the system with split preconditioner
\begin{equation}
    U^\top A U y = U^\top b, \label{eq:SplitPrecSys}
\end{equation}
whose matrix is also SPD.
System~\eqref{eq:SplitPrecSys} can then be solved with CG. The latter updates estimate $y_{\ell}$ that can be used to recover $\hat{x}_\ell := U y_\ell$.
\Cref{Algo:PCG}, the preconditioned conjugate gradients method, is equivalent to the procedure just described, but only involves solves with \(F\) and does not assume knowledge of \(U\)~\cite[p.532]{GoVa13}.
PCG updates $\hat{x}_\ell$ directly.

PCG looks for an approximate solution in the Krylov subspace 
$$
x_0 + U\mathcal{K}_{\ell}(U^\top A U , U^\top r_0),
$$
and as in~\eqref{eq:MiniProcess}, it minimizes the energy norm, 
\begin{equation}
\lVert x^\ast - \hat{x}_\ell \rVert_{A} = \min_{q \in \mathbb{P}_\ell(0)} \lVert U q (U^\top A U) U^{-1} \left( x^\ast - x_0\right) \rVert_{A}.\label{eq:error_energy}    
\end{equation}

Although there is no general method for building a good preconditioner~\cite{Preconditionning1, wathen2015preconditioning}, leveraging the convergence properties of CG on~\eqref{eq:error_energy} often leads to the following criteria: (i) $F$ should approximate the inverse of $A$, (ii) $F$ should be cheap to apply, (iii) $\kappa(U^\top A U)$ should be smaller than $ \kappa(A)$, and (iv) $U^\top A U$ should have a more favorable distribution of eigenvalues than \(A\). Note that, all four criteria only focus on $A$ and overlook the significance of the initial guess. 

\subsection{Preconditioning for a sequence of linear systems}
\label{sec:SeqOfLinSys}
In the context of~\eqref{eq:spd-sequence}, it is common to use a first level preconditioner, $F^{(1)}$, for the initial linear system, $A^{(1)} x^{(1)} = b^{(1)}$. The selection of the first-level preconditioner depends on the problem and may take into account both the physics of the problem and the algebraic structure of  $A^{(1)}$~\cite{Preconditionning1, wathen2015preconditioning, PearsonPestana2020}. To further accelerate convergence of an iterative method such as PCG on subsequent linear systems $A^{(j+1)} x^{(j+1)} = b^{(j+1)}$, one can perform a low-rank update of the most-recent preconditioner, $F^{(j)}$, leveraging information obtained from solving $A^{(j)} x^{(j)} = b^{(j)}$~\cite{AutoPrec, LMP}.

One common choice of low-rank update is to use the (approximate) spectrum of $A^{(j)}$~\cite{DataAssimilaton, SenSpectLMP, LMP}. The main idea is to capture the eigenvalues not captured by the first-level preconditioner, and cluster them to a positive quantity, typically around $1$. 

In this paper, we will consider the case where only the right-hand side is changing over the sequence of the linear systems, i.e., $A^{(j)} = A$ for all $j$. Perturbation analysis with respect to $A$ will be presented in a forthcoming paper.

\section{A scaled spectral preconditioner}\label{sec:A scaled spectral preconditioner}

We focus on the scaled spectral preconditioner, known in the literature as the deflating preconditioner~\cite{SenSpectLMP} or spectral Limited Memory Preconditioner (LMP)~\cite{LMP}, which is defined using a scaling parameter that determines the positioning of the cluster. We will provide several strategies for the choice of the scaling parameter, which has a significant impact on the convergence of PCG. 

Let us assume that $k$ largest eigenvalues of $A$, i.e. $\{\lambda_i\}_{i=1}^{k}$, are available. We define the spectral preconditioner
\begin{equation}
\label{eq:ftheta}
F_{\theta} := I_n + \sum_{i=1}^{k} \left(\frac{\theta}{\lambda_i}- 1\right) s_i s_i^\top = I_n + S_k (\theta \Lambda_k^{-1} - I_k)S_k^\top = S \begin{bmatrix}
            \theta\Lambda ^{-1}_k &  \\
            &  I_{n-k} 
       \end{bmatrix} S^\top,
\end{equation}
where $S_k := \begin{bmatrix} s_1 & \cdots & s_k \end{bmatrix}$ and $\Lambda_k := \text{diag}(\lambda_1, \ldots, \lambda_k)$.
The factor of $F_{\theta} = U_{\theta}^2$ is
\begin{equation}\label{eq:FactorizedLmp}
U_{\theta} = U_{\theta}^\top := I_n + \sum_{i=1}^{k} \left( \sqrt{\frac{\theta}{\lambda_i}}- 1\right) s_i s_i^\top 
= S 
        \begin{bmatrix}
            \sqrt{\theta}\Lambda ^{-\frac{1}{2}}_k &  \\
            &  I_{n-k} 
       \end{bmatrix} S^\top. 
\end{equation}
Preconditioner \(F_\theta\) clusters $\lambda_1, \ldots, \lambda_k$ around $\theta$, and leaves the rest of the spectrum untouched, i.e.,
\begin{equation}
    U_{\theta} A U_{\theta} = 
S \begin{bmatrix}
    \theta I_k &  \\
     &\bar{\Lambda}_{k} \\
\end{bmatrix} S^\top = \theta S_k S_k^\top + \bar{S}_{k} \bar{\Lambda}_{k} \bar{S}_{k}^\top,\label{eq:PrecMatrix}
\end{equation}
where $\bar{S}_{k} := \begin{bmatrix} s_{k+1} & \cdots & s_{n} \end{bmatrix}$ and $\bar{\Lambda}_{k} := \text{diag}(\lambda_{k+1}, \ldots, \lambda_{n})$.
As in~\eqref{eq:error_energy}, PCG minimizes
\begin{align}
\nonumber
    \lVert x^\ast - \hat{x}_\ell(\theta) \rVert_{A} & = \min_{q \in \mathbb{P}_\ell(0)} \lVert U_{\theta} q \left(U_{\theta} A  U_{\theta} \right)  U_{\theta}^{-1} \left( x^\ast - x_0\right) \rVert_{A} \\
    & = \min_{q \in \mathbb{P}_\ell(0)} \lVert q \left(U_{\theta} A  U_{\theta} \right) \left( x^\ast - x_0\right) \rVert_{A},
    \label{Eq:energy_norm_precond_system}
\end{align}
where we used $U_{\theta} q \left(U_{\theta} A  U_{\theta} \right)  U_{\theta}^{-1} = U_{\theta} U_{\theta}^{-1} q \left( U_{\theta} A  U_{\theta}\right) = q \left( U_{\theta} A  U_{\theta}\right)$. 
Using~\eqref{eq:MiniProcess} in the context of~\cref{Eq:energy_norm_precond_system}, we obtain the following result.

\begin{theorem}\label{theorem:ErrorPrecSys}
Let \(\hat{x}_\ell(\theta)\) be generated at iteration \(\ell\) of \Cref{Algo:PCG} applied to \(A x = b\) with preconditioner~\eqref{eq:ftheta}. Then, 
    \begin{equation}\label{eq:ErrorPrecond}
    \begin{aligned}
        \lVert x^\ast - \hat{x}_\ell(\theta) \rVert_{A}^2 & =
         \min_{q\in \mathbb{P}_{\ell}(0)}  \,
         \sum_{i=1}^{k} \dfrac{\eta_i^2}{\lambda_i}q(\theta)^2 + \sum_{i=k+1}^{n} \dfrac{\eta_i^2}{\lambda_i}q(\lambda_i)^2,
    \end{aligned}    
    \end{equation}
    where $\eta_i = s_i^\top r_0$ is the $i$-th component of the initial residual in the basis $S$.  
\end{theorem}
\begin{proof}
Given~\eqref{eq:PrecMatrix}, we have for any polynomial $q$,
\[
    q \left( U_{\theta} A  U_{\theta}\right) =S q \left( \begin{bmatrix}
    \theta I_k &  \\
     &\bar{\Lambda}_{k} \\
    \end{bmatrix}\right) S^\top.
\]
Since $x^\ast - x_0 =  A^{-1}r_0 = S \Lambda^{-1} S^\top r_0$,
\begin{equation}%
\label{Eq:quthetaA}
q \left( U_{\theta} A  U_{\theta}\right) (x^\ast - x_0) = S q \left( \begin{bmatrix}
    \theta I_k &  \\
     &\bar{\Lambda}_{k} \\
    \end{bmatrix}\right) \Lambda^{-1} S^\top r_0. 
\end{equation}
Substituting~\cref{Eq:quthetaA} into~\cref{Eq:energy_norm_precond_system}, we obtain the result.
\end{proof}

The scaled LMP~\cref{eq:FactorizedLmp} is typically used with $\theta = 1$. This choice is operational in numerical weather forecast~\cite{DataAssimilaton, tshimanga2007}. In the next subsections, we explore various choices for $\theta$ aiming to improve convergence properties of PCG.

\subsection{On the choice  of the scaling parameter}

The scaling parameter $\theta$, which defines the position of the cluster, is often set to $1$~\cite{DataAssimilaton, SenSpectLMP, LMP}. This choice is motivated by several factors, such as the eigenvalue distribution of $A$, the behavior of the first-level preconditioner, and the convergence behavior of PCG. 

We investigate clustering the eigenvalues at a general $\theta > 0$, which, compared with the conventional choice of $1$, results in enhanced convergence of PCG. It is important to note that the notion of ``better convergence'' may vary across different applications. For instance, in some applications, one may require high accuracy, in which case, a better convergence may be defined as a lower number of iterations. In other applications, we may want to get an approximate solution quickly, which requires to improve the convergence especially in the early iterations.
In this case, there is no guarantee that the early preconditioned iterates will provide a better reduction in the energy norm compared to the unpreconditioned iterates (\Cref{SubSec:Case1}).
For certain applications, such as numerical weather forecast, where PCG is stopped before reaching convergence due to computational budget, early convergence properties could be of critical importance. As a first direction, we will focus on the following question:
\begin{quote} 
\textit{ Is there $\theta > 0$ such that for any $x_0$,
\begin{equation}%
\label{theta_case1}
    \left\| x^{\ast}- \hat{x}_{\ell}(\theta)\right\|_{A} \leq  \left\| x^{\ast}- x_{\ell}\right\|_{A} ,\quad \ell = 1, \ldots, n?
\end{equation}}
\end{quote}
To accelerate early convergence, we will  investigate optimal choices for $\theta$ with respect to the error in the energy norm at the first iteration of PCG, i.e.,
\begin{equation*}
\min_{\theta} \Phi(\theta):= \| x^{\ast} - \hat{x}_{1}(\theta) \|_{A}^2.
\end{equation*}
We focus solely on the first iteration as it allows us to derive the optimal value of $\theta$ in closed form. 

On the other hand, for PCG, it is well known that removing eigenvalues causing convergence delay can improve the convergence rate significantly~\cite{DataAssimilaton,LMP}. This can be done by using deflation techniques, in which the aim is to ``hide" (problematic)
parts of the spectrum of  $A$ from PCG, so that the convergence rate of PCG is improved~\cite{Kahl2017,DeflatedCG}. 
Finally, our focus will be also on answering the question
\begin{quote} 
\textit{ Can we choose $\theta > 0$ such that for any $x_0$, PCG generates iterates close to those of deflation techniques?
}
\end{quote}

\subsection{\boldmath{\boldmath{\texorpdfstring{$\theta$}{theta}}} providing lower error in energy norm}
\label{SubSec:Case1}

In general, although scaled spectral preconditioning is expected to help reduce the number of iterations required to achieve convergence, \cref{theta_case1} may not hold for any choice of $\theta>0$ and all iterations $\ell$ as given by the following proposition.

\begin{proposition}\label{prop:1}
    Let $x_{1}$ be the first CG iterate when solving $Ax = b$.
    Let \(\hat{x}_1(\theta)\) be generated at the first iteration of \Cref{Algo:PCG} applied to \(A x = b\) with preconditioner~\eqref{eq:ftheta}.
    Let $x_0$ be such that $\eta_i^2 = \lambda_i$ for $i = k, k+1$ and $\eta_i = 0$ otherwise.
    Then,
    \[
    \lVert x^\ast - \hat{x}_{1}(\theta) \rVert_{A}^2 \leq \lVert x^\ast - x_{1} \rVert_{A}^2 
    \iff
    \dfrac{\lambda_{k+1}^2}{\lambda_k} \leq \theta \leq \lambda_k.
    \] 
\end{proposition}
\begin{proof}
For $\ell = 1$,~\eqref{eq:ErrorRitz} yields $\lVert x^\ast - x_1 \rVert_{A}^2  = p_1^\ast (\lambda_k)^2 + p_1^\ast(\lambda_{k+1})^2$, where
\begin{align*}
    p_1^\ast (\lambda) = 1 - \dfrac{r_0^\top r_0}{r_0^\top  A r_0} \lambda = 1 - \dfrac{\lambda_k + \lambda_{k+1}}{\lambda_k^2 + \lambda_{k+1}^2} \lambda.
\end{align*}
Similarly,~\eqref{eq:ErrorPrecond} gives $\lVert x^\ast - \hat{x}_1(\theta) \rVert_{A}^2  = q_{1, \theta}^\ast (\theta)^2 + q_{1, \theta}^\ast(\lambda_{k+1})^2$, where
\begin{align*}
    q_{1, \theta}^\ast (\lambda)  = 1 - \dfrac{r_0^\top F_{\theta} r_0}{r_0^\top F_{\theta} A F_{\theta} r_0} \lambda = 1 - \dfrac{\theta + \lambda_{k+1}}{\theta^2 + \lambda_{k+1}^2} \lambda
\end{align*}
is the polynomial that realizes the minimum. Using these relations, we obtain 
\[
\lVert x^\ast - x_1 \rVert_{A}^2 = \left(1 - \dfrac{\lambda_k + \lambda_{k+1}}{\lambda_k^2 + \lambda_{k+1}^2} \lambda_k \right)^2 + \left(1 - \dfrac{\lambda_k + \lambda_{k+1}}{\lambda_k^2 + \lambda_{k+1}^2} \lambda_{k+1} \right)^2 = \dfrac{\left(\lambda_k - \lambda_{k+1}\right)^2}{\lambda_k^2 + \lambda_{k+1}^2}
\]
and
\[
\lVert x^\ast - \hat{x}_1(\theta) \rVert_{A}^2 = \left(1 - \dfrac{\theta + \lambda_{k+1}}{\theta^2 + \lambda_{k+1}^2} \theta \right)^2 + \left(1 - \dfrac{\theta + \lambda_{k+1}}{\theta^2 + \lambda_{k+1}^2} \lambda_{k+1} \right)^2 = \dfrac{\left(\theta - \lambda_{k+1}\right)^2}{\theta^2 + \lambda_{k+1}^2}.
\]
Hence,
\begin{align*}
\dfrac{\left(\theta - \lambda_{k+1}\right)^2}{\theta^2 + \lambda_{k+1}^2} \leq \dfrac{\left(\lambda_k - \lambda_{k+1}\right)^2}{\lambda_k^2 + \lambda_{k+1}^2} \iff \dfrac{\lambda_{k+1}^2}{\lambda_k} \leq \theta \leq \lambda_k.
\end{align*}
\end{proof}

\Cref{prop:1} shows that~\cref{theta_case1} is not satisfied for all $\theta>0$.
If $\theta > 0$ lies outside of $[\lambda_{k+1}^2 / \lambda_k, \, \lambda_k]$, then \(\Vert x^\ast - \hat{x}_1(\theta) \Vert_{A} > \Vert x^\ast - x_1 \Vert_{A}\) for $x_0$ as defined in \Cref{prop:1}. 

In what comes next, we focus on the properties of $\theta$ such that \cref{theta_case1} is guaranteed for all iterations $\ell$, and for any given $x_0$.
An intuitive approach is to identify a range of $\theta$ values where the eigenvalue ratios of the preconditioned matrix are less than or equal to those of the unpreconditioned matrix,
as noted in~\cite[Lemma 1]{PolyComp}.
The next lemma shows that this property holds for $\theta \in [\lambda_{k+1}, \lambda_{k}]$, and for such choice, there exists a polynomial that promotes favorable PCG convergence.

\begin{lemma}\label{lemma:Poly_comp}
Let $\lambda_1\ge  \lambda_2\ge  \ldots \ge  \lambda_n>0$, $\ell \in \{1,\ldots,n\}$, and $k \in \{1,\ldots,\ell\}$.
For any $\theta \in [\lambda_{k+1}, \lambda_k]$, and any polynomial $p$ of degree $\ell$ such that $p(0) = 1$ and whose roots all lie in $[\lambda_n, \lambda_1]$, there exists a polynomial $q$ of degree $\ell$ such that $q(0) = 1$ and
\begin{align*}
    \lvert q(\theta) \rvert & \leq \lvert p(\lambda_i) \rvert, \quad i = 1, \ldots, k\ \\
    \lvert q(\lambda_i) \rvert & \leq \lvert p(\lambda_i) \rvert, \quad i = k+1, \ldots, n.
\end{align*}
\end{lemma}
\begin{proof}
    Let us denote $(\mu_j)_{1 \leq j \leq \ell}$ the roots of the polynomial $p$ given in decreasing order, so $p(\lambda) = \prod_{i = 1}^{\ell}\left( 1 - \frac{\lambda}{\mu_i} \right)$ for any $\lambda \geq 0$. Then, three cases may occur:

        \underline{Case 1:} For all $j \in \{1, \ldots, \ell\}$, $\mu_j < \theta$, we choose $q(\lambda) = p(\lambda)$, then simply we have for $i \in \{k+1, \ldots, n\}$, $\lvert q(\lambda_i) \rvert = \lvert p(\lambda_i) \rvert$. For $i \in \{1, \ldots, k\}$, using the property that $\mu_j < \theta \leq \lambda_i$, we obtain 
        $$1 - \frac{\lambda_i}{\mu_j} \leq 1 - \frac{\theta}{\mu_j} \leq 0.$$ 
        Thus, we have $\lvert 1 - \frac{\theta}{\mu_j} \rvert \leq \lvert 1 - \frac{\lambda_i}{\mu_j} \rvert$, and consequently $\lvert q(\theta) \rvert \leq \lvert p(\lambda_i) \rvert$.

           \underline{Case 2:} If for all $j \in \{1, \ldots, \ell\}$, $\theta \leq \mu_j$, we choose $q(\lambda) = \prod_{j = 1}^{\ell}\left( 1 - \frac{\lambda}{\theta} \right)= \left( 1 - \frac{\lambda}{\theta} \right)^l$. Then simply for $i \in \{1, \ldots, k\}$, $\lvert q(\theta) \rvert = 0 \leq \lvert p(\lambda_i) \rvert$. For $i \in \{k+1, \ldots, n\}$, using the property $\lambda_{k+1} \leq \theta \leq \mu_j$, we obtain        $$0 \leq 1 - \frac{\lambda_i}{\lambda_{k+1}} \leq 1 - \frac{\lambda_i}{\theta} \leq 1 - \frac{\lambda_i}{\mu_j}.$$        Therefore, for $i = k+1, \ldots, n$, $\lvert q(\lambda_i) \rvert \leq \lvert p(\lambda_i) \rvert$.
        
       \underline{Case 3:}  let $s \in \{1, \ldots, \ell-1\}$ such that for $j = 1, \ldots, s$, $\theta \leq \mu_j \leq \lambda_1$, and for $j = s+1, \ldots,\ell$, $\lambda_n \leq \mu_j < \theta$. Let's denote 
        $$q(\lambda) = \prod_{j = 1}^{s}\left( 1 - \frac{\lambda}{\theta} \right) \prod_{j = s+1}^{\ell}\left( 1 - \frac{\lambda}{\mu_j} \right)= \left( 1 - \frac{\lambda}{\theta} \right)^s \prod_{j = s+1}^{\ell}\left( 1 - \frac{\lambda}{\mu_j} \right) .$$
        We have $q (\theta) = 0$, so $\lvert q(\theta) \rvert \leq \lvert p(\lambda_i) \rvert$ for $i \in \{1, \ldots, k\}$. For $i \in \{k+1, \ldots, n\}$ and  $j \in \{1, \ldots, s\}$, we have 
        $$0 \leq 1 - \frac{\lambda_i}{\lambda_{k+1}} \leq 1 - \frac{\lambda_i}{\theta} \leq 1 - \frac{\lambda_i}{\mu_j},$$
        because $\lambda_{k+1} \leq \theta \leq \mu_j$. Therefore, for $i = k+1, \ldots, n$, $\lvert q(\lambda_i) \rvert \leq \lvert p(\lambda_i) \rvert$.
\end{proof}
Now, we can present a result that enables comparing the error in energy norm between the preconditioned system given by~\eqref{eq:SplitPrecSys} and the unpreconditioned system, $Ax = b$.
\begin{theorem}
\label{theorem : error_UAU_error_A}
Let $(x_{\ell})_{\ell \in  \{1,\ldots,n\}}$ and $\hat{x}_{\ell}(\theta)_{\ell \in  \{1,\ldots,n\}}$ be the sequences generated by CG and PCG with $F_{\theta}$ with $\theta \in [\lambda_{k+1}, \lambda_k]$, respectively, when solving $Ax = b$.
Assume that $\hat{x}_{0}(\theta) = x_0$. Then, for all $\ell = 1,\ldots,n$, \(\| x^\ast - \hat{x}_{\ell}(\theta) \|_{A} \leq \| x^\ast - x_{\ell} \|_{A}\).
\end{theorem}
\begin{proof}
Let $\ell \in \{1,\ldots,n\}$. From~\eqref{eq:ErrorRitz},
\begin{equation}
  \left\| x^\ast - x_{\ell}\right\|_{A}^2 = \min_{p \in \mathbb{P}_\ell(0)} \lVert p_\ell \left(A \right) \left( x^\ast - x_0\right) \rVert_{A}^2 = \sum_{i= 1}^{n}\dfrac{\eta_i^2}{\lambda_i} p_\ell^\ast(\lambda_i)^2,\label{eq:error_A}  
\end{equation}
where $\eta_i$ represents the components of the initial residual $r_0 = b - Ax_0$ in the eigenspace of $A$.
Applying \Cref{lemma:Poly_comp} to $p_{\ell}^\ast$, there exists a polynomial $q$ of degree $\ell$ with  $q(0) = 1$ such that
\begin{align*}
\lvert q(\theta) \rvert & \leq \lvert p_{\ell}^\ast(\lambda_i) \rvert, \quad i \in \{1, \ldots, k\}\\
\lvert q(\lambda_i) \rvert & \leq \lvert p_\ell^\ast(\lambda_i) \rvert, \quad i \in \{k+1, \ldots, n\}.
\end{align*}
Applying these inequalities to~\eqref{eq:error_A} yields 
\begin{align*}
   \hspace{1cm} \left\| x^\ast- x_{\ell}\right\|_{A}^2  = \sum_{i = 1}^{n} \frac{\eta_i^2}{\lambda_i} p_{\ell}^\ast(\lambda_i)^2 &  \geq \sum_{i =1}^{k} \dfrac{\eta_i^2}{\lambda_i} q(\theta)^2 + \sum_{i =k+1}^{n} \dfrac{\eta_i^2}{\lambda_i} q(\lambda_i)^2\\
          & \hspace{-2cm} \geq \min_{q\in \mathbb{P}_{\ell}(0)} 
          \left (\sum_{i=1}^{k} \dfrac{\eta_i^2}{\lambda_i}q(\theta)^2 + \sum_{i=k+1}^{n} \dfrac{\eta_i^2}{\lambda_i}q(\lambda_i)^2\right ) = \left\| x^\ast - \hat{x}_{\ell}(\theta)\right\|_{A}^2 .
\end{align*}
\end{proof}

\Cref{theorem : error_UAU_error_A} offers a range of choices for $\theta$. Next, we discuss the practical and theoretical choices from this range. Let us remind that to construct the spectral LMP~\eqref{eq:FactorizedLmp}, we are given $k$ eigenpairs. As a result, one practical choice is $\theta = \lambda_k$.
This idea is summarized in the following corollary.

\begin{corollary}
    Let \(\theta = \lambda_k\).
    Then, $\left\| x^\ast - \hat{x}_{\ell}(\lambda_k)\right\|_{A} \leq \left\| x^\ast - x_{\ell}\right\|_{A}$
    for any $x_0$ and for all $\ell \in \{1,\ldots,n\}$.    
\end{corollary}

The next theorem shows that increasing $k$ results in improved convergence.

\begin{theorem}
\label{theorem : augment_k}
Let $ 1 < k_1  \leq k_2 < n$ and $\theta_{k_1}\in [\lambda_{k_1+1}, \lambda_{k_1}]$, $\theta_{k_2}\in [\lambda_{k_2+1}, \lambda_{k_2}]$ with, $\theta_{k_2} \leq \theta_{k_1}$. 
Let $(\hat{x}_{\ell}(\theta_{k_1}))_{\ell \in  \{1,\ldots,n\}}$, $(\hat{x}_{\ell}(\theta_{k_2}))_{\ell \in  \{1,\ldots,n\}}$ be the sequences obtained from PCG iterates when solving $Ax = b$
using $F_{\theta_{k_1}}$ and $F_{\theta_{k_2}}$ respectively with an arbitrary initial guess $x_0$. 
Then, for all $\ell \in \{1,\ldots,n\}$, one has:
\begin{equation*}
\left\| x^\ast- \hat{x}_{\ell}(\theta_{k_2})\right\|_{A} \leq  \left\| x^\ast - \hat{x}_{\ell}(\theta_{k_1})\right\|_{A}.
\end{equation*}
\end{theorem}
\begin{proof}
    The eigenvalues of the preconditioned matrix using $F_{\theta_{k_1}}$  and $F_{\theta_{k_2}}$
    are given in decreasing order respectively as  
\[
       \rho_i = \begin{cases}
        \theta_{k_1} &   i \in \{1, \ldots, k_1\}\\
        \lambda_i & \text{otherwise},
    \end{cases}
   \quad \text{and} \quad
     \widetilde{\rho}_i = \begin{cases}
        \theta_{k_2} &   i \in \{1, \ldots, k_2\}\\
        \lambda_i & \text{otherwise}.
    \end{cases}
\]
As $k_1 < k_2$, it follows that  $\widetilde{\rho}_{k_2} \leq \rho_{k_1} = \theta_{k_1}$. Therefore, $\widetilde{\rho}_{i}$ can be expressed as a function of $\rho_i$ as 
\[
       \widetilde{\rho}_i = \begin{cases}
        \theta_{k_2} \in [\rho_{k_2 +1},\rho_{k_2}]&   i \in \{1, \ldots, k_2\}\\
        \rho_i & \text{otherwise}.
    \end{cases}
\]
Using \Cref{lemma:Poly_comp}, for the polynomial $q_{\ell, \theta_{k_1}}^{\ast}$, there exists a polynomial $q$ of degree $\ell$ with  $q(0) = 1$, such that for $i \in \{1, \ldots,n\}$,
\begin{align*}
\lvert q(\theta_{k_2}) \rvert & \leq \lvert q_{\ell, \theta_{k_1}}^\ast(\rho_i) \rvert, \quad   i \in \{1, \ldots, k_2\}\\
\lvert q(\rho_i) \rvert & \leq \lvert q_{\ell, \theta_{k_1}}^\ast(\rho_i) \rvert, \quad  i \in \{k_2+1, \ldots, n\}\\
\end{align*}
Applying this result to~\eqref{eq:ErrorPrecond} yields that
\begin{align*}
   \hspace{1cm} \left\| x^\ast- \hat{x}_{\ell}(\theta_{k_1})\right\|_{A}^2 & = \sum_{i = 1}^{n} \frac{\eta_i^2}{\lambda_i} q_{\ell, \theta_{k_1}}^\ast(\rho_i)^2 \\
          & \hspace{-1cm} \geq \sum_{i =1}^{k_2} \dfrac{\eta_i^2}{\lambda_i} q(\theta_{k_2})^2 + \sum_{i =k_2+1}^{n} \dfrac{\eta_i^2}{\lambda_i} q(\rho_i)^2\\
          &  \hspace{-1cm} \geq \min_{q\in \mathbb{P}_{\ell}(0)} 
          \left (\sum_{i=1}^{k_2} \dfrac{\eta_i^2}{\lambda_i}q(\theta_{k_2})^2 + \sum_{i=k_2+1}^{n} \dfrac{\eta_i^2}{\lambda_i}q(\lambda_i)^2\right )= \left\| x^\ast - \hat{x}_{\ell}(\theta_{k_2})\right\|_{A}^2 .
\end{align*}
\end{proof}
One can see that $k_1 < k_2 \ \Longrightarrow \ \theta_{k_2} \leq \theta_{k_1}$, since $\lambda_i$ are in decreasing order. In addition, when $k_1 = k_2$, \Cref{theorem : augment_k} shows that $\lambda_{k_1+1}$ is the best choice in $[\lambda_{k_1+1}, \lambda_{k_1}]$ in terms of reducing the error with respect to the unpreconditioned system.

\subsection{Optimal choice for \boldmath{\texorpdfstring{$\theta$}{theta}} with respect to the initial residual}\label{SubSec:Case2}
Our objective is to determine the value of $\theta$ that minimizes the energy norm of the error at the initial iterate. This will provide us with the optimal reduction at the first iterate,
\begin{equation}
    \thetares \in \arg\min_{\theta > 0} \Phi(\theta):= \left\| x^\ast - \hat{x}_{1}(\theta)\right\|_{A}^2.\label{eq:min_def}
\end{equation}

The expression for $\thetares$ is stated in the following theorem.

\begin{theorem}
Let $r_0 = b - Ax_0$.
The unique $\lambda_n \leq \thetares \leq \lambda_{k+1}$ satisfying~\eqref{eq:min_def} is
\begin{equation}
\label{eq:theta_r}
    \thetares := \dfrac{\sum_{i=k+1}^{n}\lambda_i\eta_i^2}{\sum_{i=k+1}^{n}\eta_i^2} = \dfrac{r_0^\top A r_0 - r_0 S_k \Lambda_k S_k^\top r_0}{r_0^\top r_0 - r_0^\top S_k S_k^\top r_0}.
\end{equation}
\end{theorem}

\begin{proof}
First, \Cref{theorem:ErrorPrecSys} implies
\begin{equation} \left\| x^\ast - \hat{x}_{1}(\theta)\right\|_{A}^2 = \sum_{i=1}^{k}\dfrac{\eta_i^2}{\lambda_i} q_{1, \theta}^\ast (\theta)^2 + \sum_{i=k+1}^{n}\dfrac{\eta_i^2}{\lambda_i} q_{1, \theta}^\ast(\lambda_i)^2\label{eq:phi} \end{equation}
 where $\eta_i = s_i^\top r_0$ and  $
q_{1, \theta}^\ast (\lambda) = 1 - \dfrac{r_0^\top F_{\theta} r_0}{r_0^\top F_{\theta} A F_{\theta} r_0} \lambda$. 
Using~\eqref{eq:ftheta}, we obtain 
\begin{equation}
    r_0^\top F_{\theta} r_0 = \theta \sum_{i=1}^{k}\frac{\eta_i^2}{\lambda_i} +  \sum_{i=k+1}^{n}\eta_i^2 \quad \mbox{and} \quad   r_0^\top F_{\theta} A F_{\theta} r_0 = \theta^2 \sum_{i=1}^{k}\frac{\eta_i^2}{\lambda_i} +  \sum_{i=k+1}^{n}\lambda_i \eta_i^2.
    \label{eq:rFr}
\end{equation}
Then, for all $\theta > 0$, $\Phi(\theta)$ simplifies to
\begin{equation*}
 \Phi(\theta)=
a_1 \left(\dfrac{a_2 \theta -a_3}{a_1 \theta^2 + a_3}\right)^2+\sum_{i=k+1}^n \dfrac{\eta_i^2}{\lambda_i} \left(1 - \dfrac{a_1 \theta +  a_2}{a_1 \theta^2 + a_3} \lambda_i\right)^2,
\end{equation*}
where $a_1 = \displaystyle\sum_{i=1}^{k}\dfrac{\eta_i^2}{\lambda_i}$, $a_2 = \displaystyle\sum_{i=k+1}^{n}\eta_i^2$ and $a_3 = \displaystyle \sum_{i=k+1}^{n}\lambda_i \eta_i^2$.  The derivative of $\Phi$ is 
\begin{equation*}
        \Phi'(\theta) 
        = \dfrac{2 a_1}{(a_1\theta^2+a_3)^3} (a_2\theta - a_3) (a^2 \theta^3 + a_1a_2 \theta^2 +a_1a_3 \theta + a_2a_3).
\end{equation*}
Since $\Phi'(\theta) < 0$ on $]0, \frac{a_3}{a_2}[$ and $\Phi'(\theta)  > 0$ on $]\frac{a_3}{a_2}, +\infty[$, then $\frac{a_3}{a_2}$ is the global minimizer of $\Phi$ on $\R^{*}_+$ and is unique. Hence,
\[\thetares = \arg\min_{\theta>0} \Phi(\theta)= \dfrac{a_3}{a_2} = \dfrac{\sum_{i=k+1}^{n}\lambda_i\eta_i^2}{\sum_{i=k+1}^{n}\eta_i^2}.
\]
Moreover,
\[\lambda_n=\dfrac{\sum_{i=k+1}^{n}\lambda_n\eta_i^2}{\sum_{i=k+1}^{n}\eta_i^2}\leq \thetares  \leq \dfrac{\sum_{i=k+1}^{n}\lambda_i\eta_i^2}{\sum_{i=k+1}^{n}\eta_i^2} = \lambda_{k+1}.
\]
The expression for $\thetares$ can be rewritten in terms of $S_k$, $\Lambda_k$, and $r_0$ as follows:
\begin{align*}
        \thetares =\dfrac{\sum_{i=1}^{n}\lambda_i\eta_i^2 - \sum_{i=1}^{k}\lambda_i\eta_i^2} {\sum_{i=1}^{n}\eta_i^2 - \sum_{i=1}^{k}\eta_i^2} = \dfrac{r_0^\top A r_0 - r_0 S_k \Lambda_k S_k^\top r_0}{r_0^\top r_0 - r_0^\top S_k S_k^\top r_0}.
\end{align*}
\end{proof}

Note that $\thetares$ can be interpreted as the center of mass for the remaining part of the spectrum in which the weights are determined by $\eta_i^2$, i.e.
 $$
 \sum_{i=k+1}^{n} \eta_i^2 (\thetares - \lambda_i ) =0.
 $$ 

Let us now look at the first iterate,
\begin{equation}
    \hat{x}_1(\thetares) = x_0 + \frac{r_0^\top F_{\thetares} r_0}{r_0^\top F_{\thetares} A F_{\thetares} r_0} F_{\thetares} r_0,\label{eq:x1}
\end{equation}
to better understand the effect of 
$\thetares$.
Using~\eqref{eq:rFr} and the value of $\thetares$,

\[
\frac{r_0^\top F_{\thetares} r_0}{r_0^\top F_{\thetares} A F_{\thetares} r_0} = \frac{\sum_{i=k+1}^{n}\eta_i^2}{\sum_{i=k+1}^{n}\lambda_i\eta_i^2} = \dfrac{1}{\thetares}.
\]
Therefore,~\eqref{eq:x1} simplifies to
\begin{equation*}
    \hat{x}_1(\thetares) = x_0 + \frac{1}{\thetares} \left(\bar{S}_{k}\bar{S}_{k}^\top + \thetares S_{k} \Lambda_k^{-1} S_{k}^\top \right) r_0= x_0 + S_{k} \Lambda_k^{-1} S_{k}^\top r_0 + \frac{1}{\thetares}\bar{S}_{k}\bar{S}_{k}^\top r_0.
\end{equation*}
Then, the residual of the first iteration 
is given by
\begin{equation}\label{eq:r_1_mass}
    \begin{aligned}
    b - A \hat{x}_1(\thetares) &= r_0 - S_k S_k^\top r_0 - \frac{1}{\thetares} \bar{S}_{k} \bar{\Lambda}_{k} \bar{S}_{k}^\top r_0 = \bar{S}_{k} \bar{S}_{k}^\top r_0  - \frac{1}{\thetares} \bar{S}_{k} \bar{\Lambda}_{k} \bar{S}_{k}^\top r_0. 
\end{aligned}
\end{equation}
Given~\eqref{eq:r_1_mass}, we conclude that, from the first iteration, we can remove all components of the residual with respect to $S_k$, see \Cref{appendix:DefCG}. We now provide an upper bound for the error in the energy norm for later iterations, $\ell > 1$, beginning with $\hat{x}_1(\thetares)$. With this initial point, we ensure that all iterates yield a residual within $\text{Span}(\bar{S}_{k})$.

\begin{theorem}
    Let $\hat{x}_{\ell}(\thetares)$ be the $\ell$-th iterate obtained from PCG when solving $Ax = b$ using the preconditioner $F_{\thetares}$ with an arbitrary initial guess $x_0$. Let $x_{\ell}^{\text{Init}}$ be the $\ell$-th iterate generated by CG for solving $A x = b$ starting from $\hat{x}_1(\thetares)$ as defined in~\eqref{eq:x1}.
Then, for all $\ell \in \{1, \ldots, n\}$, $\left\| x^\ast- \hat{x}_{\ell + 1}(\thetares)\right\|_{A} \leq  \left\| x^\ast- x_{\ell}^{\text{Init}}\right\|_{A}$.
\end{theorem}

\begin{proof}
    From~\eqref{eq:r_1_mass}, the components of $b - A \hat{x}_1(\thetares)$ in the eigenspace of $A$ are
     \begin{equation*}
        0 \quad (i = 1, \ldots, k), \quad \textup{and} \quad
        \eta_i ( 1 -  \lambda_i / \thetares) \quad (i > k).\\
     \end{equation*}
Thus, 
    \begin{equation}
        \left\| x- x_{\ell}^{\text{Init}}\right\|_{A}^2 = \sum_{i = k+1}^{n} \dfrac{\eta_i^2}{\lambda_i}  \left( 1 -  \dfrac{\lambda_i}{\thetares}\right)^2 p_\ell^{\ast, {\text{Init}}}(\lambda_i)^2,
    \end{equation}
where $p_\ell^{\ast, {\text{Init}}}$ is the polynomial that minimizes $p \mapsto \lVert p \left(A \right) \left( x^\ast - \hat{x}_1(\thetares)\right) \rVert_{A}^2$ over $\mathbb{P}_{\ell}(0)$.\\ 
Define 
\[
\Bar{q}(\lambda) = \left( 1 -  \dfrac{\lambda}{\thetares}\right) p_\ell^{\ast, {\text{Init}}}(\lambda),
\] 
and note that $\Bar{q} \in \mathbb{P}_\ell(0)$. Now we have 
\begin{align*}
    \left\| x^\ast - \hat{x}_{\ell + 1}(\thetares)\right\|_{A}^2 &= \min_{q\in \mathbb{P}_{\ell+1}(0)} 
          \left (\sum_{i=1}^{k} \dfrac{\eta_i^2}{\lambda_i}q(\thetares)^2 + \sum_{i=k+1}^{n} \dfrac{\eta_i^2}{\lambda_i}q(\lambda_i)^2\right )\\
          & \leq \sum_{i=1}^{k} \dfrac{\eta_i^2}{\lambda_i}\Bar{q}(\thetares)^2 + \sum_{i=k+1}^{n} \dfrac{\eta_i^2}{\lambda_i}\Bar{q}(\lambda_i)^2\\
          & = \sum_{i = k+1}^{n} \dfrac{\eta_i^2}{\lambda_i}  \left( 1 -  \dfrac{\lambda_i}{\thetares}\right)^2 p_\ell^{\ast, {\text{Init}}}(\lambda_i)^2 = \left\| x- x_{\ell}^{\text{Init}}\right\|_{A}^2.
\end{align*}

\end{proof}

    Note that, one can interpret $\hat{x}_{1}(\thetares)$ as the first iteration of CG when solving the unpreconditioned system, starting from $x_0 + S_{k} \Lambda_k^{-1} S_{k}^\top r_0$, since the search direction at the first iteration is equal to:
\begin{equation}\label{eq:residual_barx0}
       b -  A \left( x_0 + S_{k} \Lambda_k^{-1} S_{k}^\top r_0 \right)
    = b - A x_0 - S_{k} S_{k}^\top r_0
    = r_0 - S_{k} S_{k}^\top r_0 
    = \bar{S}_{k} \bar{S}_{k}^\top r_0, 
    \end{equation}
    and the step-length $\alpha_0$ is given as 
    \begin{equation*}
  \alpha_0 = \dfrac{1}{\thetares} = \dfrac{r_0^\top \bar{S}_{k}\bar{S}_{k}^\top r_0}{r_0^\top \bar{S}_{k}^\top \bar{S}_{k} A \bar{S}_{k}\bar{S}_{k}^\top r_0}.
    \end{equation*}
This highlights the strong connection between preconditioning, CG with different initial point and deflation techniques~\cite{DeflatedCG,tshimanga2007}. This connection will be explored in detail in the next subsection, providing another choice for the scaling parameter.
\subsection{\boldmath{ \texorpdfstring{$\theta$}{theta}  } as the mid-range between \texorpdfstring{$\lambda_k$}{lambdak} and \texorpdfstring{$\lambda_n$}{lambdan}}
\label{SubSec:median}
We focus now on choosing a scaling parameter $\theta$ to obtain approximate iterates to those of deflated CG (see \Cref{al:Deflated-CG}). The deflation technique, with $S_k$ as the deflation subspace, is equivalent to standard CG applied to $A x = b$ with initial guess
\[
x_0^{\text{Def}} = x_0 + S_k \Lambda_k^{-1} S_k^\top (b - A x_0).
\]
From~\eqref{eq:residual_barx0}, the residual of $x_0^{\text{Def}}$ is given as
\begin{equation*}
   b - A x_0^{\text{Def}} = \bar{S}_{k} \bar{S}_{k}^\top r_0.
\end{equation*}
One can see that this initial guess gives a residual which is an orthogonal projection of $r_0$ onto $\text{span}(\bar{S}_{k})$, so that the $\ell$-th iterate of CG, $x_\ell^{\text{Def}}$, starting with $x_0^{\text{Def}}$ satisfies
\begin{equation*}
    \left\| x^\ast - x_\ell^{\text{Def}}\right\|_{A}^2 =
    \min_{q\in \mathbb{P}_{\ell}(0)} 
    \left (\sum_{i=k + 1}^{n} \dfrac{\eta_i^2}{\lambda_i}{q}(\lambda_i)^2 \right ).          
\end{equation*}
We now provide the main result of this section.
\begin{theorem}\label{theorem:midrange}
Let $\hat{x}_\ell(\theta)$ be the $\ell$-th iterate obtained from PCG iterates when solving $Ax = b$ using $F_{\theta}$ starting from an arbitrary initial guess $x_0 \in \R^n$.
Let $x_\ell^{\text{Def}}$ be the $\ell$-th iterate generated with CG when solving $Ax = b$ starting with $x_0^{\text{Def}} = x_0 + S_k \Lambda_k^{-1} S_k^\top (b - A x_0)$. Then, in exact arithmetic,
\begin{equation}\label{eq:upperboundmed}
 \left\| x^\ast - x_{\ell + 1}^{\text{Def}}\right\|_{A} \leq \left\| x^\ast - \hat{x}_{\ell + 1}(\theta)\right\|_{A} \leq  \dfrac{\alpha(\theta)}{\theta} \left\| x^\ast - x_\ell^{\text{Def}}\right\|_{A},
\end{equation}
 with 
$\alpha(\theta) = \max\left(|\lambda_{k+1} - \theta|,| \theta - \lambda_n |\right).$
\end{theorem}
 
\begin{proof}
    Let us start by showing the first inequality. From \Cref{theorem:ErrorPrecSys}
\begin{align*}
    \left\| x^\ast - \hat{x}_{\ell + 1}(\theta)\right\|_{A}^2 
    & = \sum_{i = 1}^{k} \frac{\eta_i^2}{\lambda_i} q_{\ell + 1, \theta}^\ast(\theta)^2 + \sum_{i = k + 1}^{n} \frac{\eta_i^2}{\lambda_i} q_{\ell + 1, \theta}^\ast(\lambda_i)^2\\
    & \geq \sum_{i = k + 1}^{n} \frac{\eta_i^2}{\lambda_i} q_{\ell + 1, \theta}^\ast(\lambda_i)^2\\
    & \geq \min_{q\in \mathbb{P}_{\ell + 1}(0)} 
          \left (\sum_{i=k + 1}^{n} \dfrac{\eta_i^2}{\lambda_i}{q}(\lambda_i)^2 \right ) 
          = \left\| x^\ast - x_{\ell + 1}^{\text{Def}}\right\|_{A}^2.
\end{align*}    
Now, to prove the second inequality, we consider $p_\ell^{\ast, \text{Def}}$ the polynomial that minimizes $p \mapsto \lVert p \left(A \right) \left( x^\ast - x_0^{\text{Def}}\right) \rVert_{A}^2$ over $\mathbb{P}_{\ell}(0)$., i.e.,
$$
\left\| x^\ast - x_{\ell}^{\text{Def}}\right\|_{A}^2 = \sum_{i=k + 1}^{n} \dfrac{\eta_i^2}{\lambda_i}{p_\ell^{\ast, \text{Def}}}(\lambda_i)^2. 
$$
Consider $\widetilde{q}_{\ell + 1} \in \mathbb{P}_{\ell + 1}(0)$ such as for all $\lambda \in \R$,$
\widetilde{q}_{\ell + 1}(\lambda) = \left(1 - \dfrac{\lambda}{\theta}\right) p_\ell^{\ast, \text{Def}}(\lambda)$.
Hence,
\begin{align*}
    \left\| x^\ast - \hat{x}_{\ell + 1}(\theta)\right\|_{A}^2 
    & = \sum_{i = 1}^{k} \frac{\eta_i^2}{\lambda_i} q_{\ell + 1, \theta}^\ast(\theta)^2 + \sum_{i = k + 1}^{n} \frac{\eta_i^2}{\lambda_i} q_{\ell + 1, \theta}^\ast(\lambda_i)^2\\
    & \leq  \sum_{i = 1}^{k} \frac{\eta_i^2}{\lambda_i} \widetilde{q}_{\ell + 1}(\theta)^2 + \sum_{i = k + 1}^{n} \frac{\eta_i^2}{\lambda_i} \widetilde{q}_{\ell + 1}(\lambda_i)^2\\
    & = \sum_{i = k + 1}^{n} \frac{\eta_i^2}{\lambda_i} p_\ell^{\text{Def}, \ast}(\lambda_i)\left(1 - \dfrac{\lambda_i}{\theta}\right)^2\\
    & \leq \max_{k+1 \leq i \leq n}\left(1 - \dfrac{\lambda_i}{\theta}\right)^2 \left\| x^\ast - x_{\ell}^{\text{Def}}\right\|_{A}^2 = \dfrac{\alpha(\theta)}{\theta} \left\|  x^\ast- x_{\ell}^{\text{Def}}\right\|_{A}^2.
\end{align*}  
\end{proof}

Choosing $\theta>0$ such that $\alpha(\theta)/\theta > 1$  in~\eqref{eq:upperboundmed} would give a pessimistic upper bound.  For a better bound, we select $\theta>0$ such that $\alpha(\theta)/\theta \leq 1$, which is equivalent to impose $\theta \ge \lambda_{k+1}/2$.
The value of $\theta$ that minimizes $\alpha(\theta)/\theta$ is  $\theta^{\ast} = (\lambda_{k+1}+\lambda_n)/2$. 

Given that $\lambda_{k+1}$ is unknown, and $\lambda_n$ can be predetermined in various applications, e.g., in data assimilation problems $\lambda_n = 1$, a practical approach for selecting $\theta$ (the closest to $\theta^{\ast}$) is by choosing the average between the $\lambda_k$ and 
 $\lambda_n$, i.e., $\thetamidrange = (\lambda_{k}+\lambda_n)/2$, for which we have $\alpha(\thetamidrange)/\thetamidrange = (\lambda_k - \lambda_n)/(\lambda_{k} + \lambda_n)< 1$.
 Note that the choice $\theta = \lambda_k$ yields in~\eqref{eq:upperboundmed} to a worst upper bound compared to $\thetamidrange$, i.e., $\alpha(\lambda_k)/\lambda_k > \alpha(\thetamidrange)/\thetamidrange$.

\subsection{Discussion}

The analysis in this section raises two key questions. The first is: why use a scaled spectral preconditioner when we know that deflated CG iterations using the deflated subspace $S_k$, or using an initial guess as defined in~\eqref{eq:x1}, produce better results in exact arithmetic (see \Cref{theorem:midrange})? The assumption in this section is that the eigenpairs used to construct the deflated subspace or the initial guess are exact, ensuring that components of the initial residual within the eigenspace of \( S_k \) are eliminated. However, when an approximate eigen-spectrum is used, such as the eigen-spectrum of \( A \) is applied to solve a system involving a perturbed matrix, \( \widetilde{A} \), the initial guess may fail to remove the components of the initial residual within the eigenspace of \( \widetilde{A} \). For instance, consider the perturbed matrix  \( \widetilde{A} = A + E \),
\( A \) is modified by a small perturbation matrix \( E \). This results in the following expression:
\[
b - \widetilde{A} x_0^{\text{Def}} = b - A x_0^{\text{Def}} + E x_0^{\text{Def}},
\]
where the value of \( b - A x_0^{\text{Def}} \) from~\eqref{eq:residual_barx0} becomes:
\(
b - \widetilde{A} x_0^{\text{Def}} = \bar{S}_{k} \bar{S}_{k}^\top \left( b - A x_0 \right) + E x_0^{\text{Def}}.
\)

This illustrates that the perturbation $E$ introduces additional components to the residual, which the initial guess fails to fully eliminate, unlike in the exact case. When the perturbation exists, we show in numerical experiments that using a scaled spectral LMP becomes advantageous over deflated CG. 

The second question is: why not combine the initial guess~\eqref{eq:x1} with the scaled spectral LMP using \( \theta = 1 \). When the initial guess fails to eliminate components of the initial residual within the eigenspace of \( \widetilde{A} \), these components influence the convergence of PCG. Their impact on the energy norm of the error can be reduced by appropriately positioning the largest eigenvalues. 


\section{Numerical Experiments}
\label{sec:NumExp}

In this section, we illustrate the performance of the scaled spectral LMP, as defined in~\eqref{eq:FactorizedLmp}, within the context of a nonlinear weighted least-squares problem arising in data assimilation, i.e.,
\begin{equation}
\label{Prob:nonlinear}
\min_{w_0 \in \R^n}  f(w_0) = \min_{w_0 \in \R^n} \frac{1}{2} \|w_0 - w_b\|_{B^{-1}}^2 + \frac{1}{2} \sum_{i=1}^{N_t} \|y_i - \mathcal{H}_i(\mathcal{M}_{t_0, t_i}(w_0))\|_{R_i^{-1}}^2.
\end{equation}
Here, $w_0 = w(t_0)$, is the state at the initial time $t_0$, for instance temperature value, $w_b \in \R^n$ is a priori information at time $t_0$ and $y_i \in \R^{m_i}$ represents the observation vector at time $t_i$ for $i=1, \ldots, N_t$. $\mathcal{M}_{t_0, t_i}(\cdot)$ is a nonlinear physical dynamical model which propagates the state $w_0$ at time $t_0$ to the the state $w_i$ at time $t_i$ by solving the partial differential equations. 
$\mathcal{H}_i(\cdot)$ maps the state vector $w_i$ to a $m_i$-dimensional vector representing the state
vector in the observation space.
$B \in \R^{n \times n}$, $R_i \in \R^{m_i \times m_i}$ are symmetric
positive definite error covariance matrices corresponding to the a priori and observation model error, respectively. 

The TGN method~\cite{gratton2007approximate} is widely used 
to solve the nonlinear optimization problem~\eqref{Prob:nonlinear}. At each iteration $j$ of the TGN method, the linearized least-squares approximation to the nonlinear least-squares problem~\eqref{Prob:nonlinear} is solved. This quadratic cost function at the $j$-th iterate is formulated as
\begin{equation} 
\label{Prob:quadratic}
      Q^{(j)}(s) =  \frac{1}{2} \left \| s - (w_b - w_0^{(j)}) \right \|_{B^{-1}}^2
      +  \frac{1}{2} \sum_{i=1}^{N_t} \|G_i^{(j)}s_i - d_i^{(j)} \|_{R_i^{-1}}^2,
\end{equation}
where $s \in \R^{n}$ , $d_i^{(j)} = y_i - \mathcal{G}_i(w_0^{(j)})$ with $\mathcal{G}_i(w_0^{(j)}) = \mathcal{H}_i(\mathcal{M}_{t_0, t_i}(w_0^{(j)}))$ and $G_i^{(j)}$ represents the Jacobian of $\mathcal{G}_i$ at a given iterate $w_0^{(j)}$. The quadratic cost function~\eqref{Prob:quadratic} is minimized with respect to $s$ which is then used to update the current iterate, i.e. $w_0^{({j+1})} = w_0^{({j})} + s^{(j)}$, where $s^{(j)}$ is an approximate solution of the problem~~\eqref{Prob:quadratic}. This process continues till the convergence criterion is met. For large scale problems with computationally expensive models $\mathcal{M}_{t_0, t_i}(\cdot)$, a limited number of TGN iterations are applied. The solution to the quadratic problem~\eqref{Prob:quadratic} can be found by solving
\begin{equation}\label{eq:linear_system}
\left( B^{-1} + (G^{(j)})^\top R^{-1} G^{(j)} \right) s = B^{-1} (w_b - w_0^{(j)}) - (G^{(j)})^\top R^{-1} d^{(j)}. 
\end{equation}
where $d^{(j)}$ is a $m$-dimensional concatenated vector of $d_i^{(j)}$ with $m = \sum_{i = 1}^{N_t} m_i$, $G^{(j)} \in \R^{m \times n}$ represents a concatenation of $G_i^{(j)} \in \R^{m_i\times n}$, and $R \in \R^{m \times m}$ is a block diagonal matrix, i.e. $R = \mbox{diag}(R_1, \ldots, R_N)$. The matrix $B^{-1} + (G^{(j)})^\top R^{-1} G^{(j)}$ is SPD, matrix-vector products with it are accessible only through operators, and $n$ can be large for data assimilation problems. Hence, CG is widely used to solve such systems.

Let us assume that a square root factorization of $B=LL^\top$ is available. The linear system~\eqref{eq:linear_system} can be then preconditioned by using this \textit{first-level} split preconditioner,
\begin{equation}
    \begin{aligned}
        &\left( I_n + L^\top (G^{(j)})^\top R^{-1} G^{(j)} L  \right) x = L^\top \left( B^{-1}(w_b - w_0^{(j)}) - (G^{(j)})^{\top} R^{-1} d^{(j)} \right).
    \end{aligned} \label{eq:linear_subproblem}
\end{equation}
CG at the $\ell$-th iteration provides an approximate solution $x_{\ell}^{(j)}$ which is then used to obtain an approximate solution of the linear system~\eqref{eq:linear_system}, i.e. $s_{\ell}^{(j)} = Lx_{\ell}^{(j)}$. In operational data assimilation problems, in general $m \ll n$. Consequently, the preconditioned matrix $A^{(j)} = I_n + L^\top (G^{(j)})^\top R^{-1} G^{(j)} L $ has $n-m$ eigenvalues clustered around $1$, while the remaining eigenvalues are greater than $1$.  

Since in the context of TGN, a sequence of closely related linear systems is solved, it is common to update the first-level preconditioner $L$ by using approximate eigenspectrum of the previous linear system~\cite{DataAssimilaton, LMP}. Let us denote $b^{(j)} := L^\top \left( B^{-1}(w_b - w_0^{(j)}) - (G^{(j)})^{\top} R^{-1} d^{(j)} \right)$. For $j=1$, CG Algorithm~\ref{Algo:CG} solves the linear system $A^{(1)} x = b^{(1)}$, for the variable $x$. Using the recurrences of CG, we can easily compute approximate eigenpairs of $A^{(1)}$ (see ~\cite[p.174]{YoussSaad} for more details). These pairs can then be used to construct a second-level preconditioner, $U_{\theta_1}^{(1)}$, by using the formula~\eqref{eq:FactorizedLmp}. Consequently, $(U_{\theta_1}^{(1)})^2$ is an approximation to the inverse of the matrix $A^{(1)}$. Then, assuming that $A^{(2)}$ is close to the matrix $A^{(1)}$, for $j=2$, 
CG Algorithm~\ref{Algo:CG} is applied to the preconditioned system, $U_{\theta}^{(1)} A^{(2)}U_{\theta_1}^{(1)} x = U_{\theta_1}^{(1)} b^{(2)}$. The approximate solution at $\ell$-iterate is obtained from the relation $s_{\ell}^{(2)} = L U_{\theta_1}^{(1)} x_{\ell}^{(2)}$. At the end of the CG, we can obtain approximate eigenpairs of $U_{\theta}^{(1)} A^{(2)}U_{\theta_1}^{(1)}$ and use it to construct a preconditioner for the next linear system. At the $j$-th outer loop of TGN, CG is applied to the preconditioned linear system:
\begin{equation}
\label{eq:precond_systems}
( U_{\theta_{j-1}}^{(j-1)} \ldots U_{\theta_1}^{(1)} A^{(j)} U_{\theta_1}^{(1)} \ldots U_{\theta_{j-1}}^{(j-1)}) \; x = U_{\theta_{j-1}}^{(j-1)} \ldots U_{\theta_1}^{(1)} b^{(j)},
\end{equation}
and the approximate solution to~\eqref{eq:linear_system} is obtained from
$s_{\ell}^{(j)} = L U_{\theta_{j-1}}^{(j-1)} \ldots U_{\theta_1}^{(1)} x_{\ell}^{(j)}$.

\subsection{Setup}
In our numerical experiments, we use the Lorenz-96~\cite{lorenz1996predictability} model as the physical dynamical system, $\mathcal{M}_{t_0, t_i}(\cdot)$, which is commonly used as a reference model in data assimilation. The observation operator $\mathcal{H}(\cdot)$ is defined as a uniform selection operator, meaning $\mathcal{H}(x)$ extracts a subset of $x$ that is uniformly selected.
$B$ is chosen as a discretized diffusion operator with a standard deviation $\sigma_b = 0.8$~\cite{goux2024impact}.  We consider $R_1 = R_{2} = \sigma_r^2 I_m$ with $\sigma_r = 0.2$.
    We choose $n = 1000$ and $N_t= 2$, and we consider two different scenarios, with a different number of observations: (1) \textit{LowObs} with $m_1 = m_2 = 150$ and (2) \textit{HighObs} with $m_1 = m_2 = 300$.
For both cases, $2$ outer loops are performed within TGN. CG is applied to the first linear system \(A^{(1)} x = b^{(1)}\) with $100$ iterations. Then, approximate largest eigen-pairs of $A^{(1)}$, $(S_k, \Lambda_k)$, are computed and selected based on convergence criteria with a tolerance of \(\varepsilon = 10^{-3}\) (See~[Section 1.3]\cite{tshimanga2007} for further details). With this criteria, the number of selected eigen-pairs is $45$ in the \textit{LowObs} case and $26$ in the \textit{HighObs} case. Using these pairs, the scaled LMP, $U_{\theta_1}^{(1)}$, is applied as a preconditioner for $j=2$. Matrix-vector products with the preconditioner are carried out via an operator using the selected pairs, meaning the preconditioner matrix is not explicitly constructed. 

\subsection{Numerical Results}
In this section, we present numerical results only for the second outer loop ($j=2$) of the TGN method. We compare the performance of the methodologies of Table~\ref{tab:methods} in terms of convergence rate and computational cost.

\begin{table}[ht]
\centering
 \resizebox{\textwidth}{!}{
\begin{tabular}{lll}
\hline
Method & Description & Initial guess \\ \hline
\textbf{BPrec} & \Cref{Algo:CG} applied to~\eqref{eq:linear_subproblem} & $x_0 = 0$ \\ 
\textbf{sLMP-Base} & \Cref{Algo:CG} applied to~\eqref{eq:precond_systems}, $\theta_1=1$ & $x_0 = 0$ \\ 
\textbf{Init-sLMP-Base} & \Cref{Algo:CG} applied to~\eqref{eq:precond_systems}, $\theta_1=1$ & $x_0 = U_{\theta_1}^{-1}S_k \Lambda_k^{-1} S_k^\top b^{(2)}$ \\ 
\textbf{sLMP-$\lambda_k$} & \Cref{Algo:CG} applied to~\eqref{eq:precond_systems}, $\theta_1=\lambda_k$ & $x_0 = 0$  \\ 
\textbf{sLMP-$\theta_r$} &
\Cref{Algo:CG} applied to~\eqref{eq:precond_systems}, $\theta_1=\theta_r$ & $x_0 = 0$  \\
\textbf{sLMP-$\theta_m$} & \Cref{Algo:CG} applied to~\eqref{eq:precond_systems}, $\theta_1=(\lambda_k + 1)/2$ & $x_0 = 0$  \\
\textbf{DefCG} & \Cref{al:Deflated-CG} applied to~\eqref{eq:linear_subproblem}, $W = S_k$ & $x_{-1} =0$ \\ \hline
\end{tabular}
}
\caption{Description of methods used in the numerical experiments}
\label{tab:methods}
\end{table}
\begin{figure}[ht]
\centering
    \includegraphics[width=.95\linewidth]{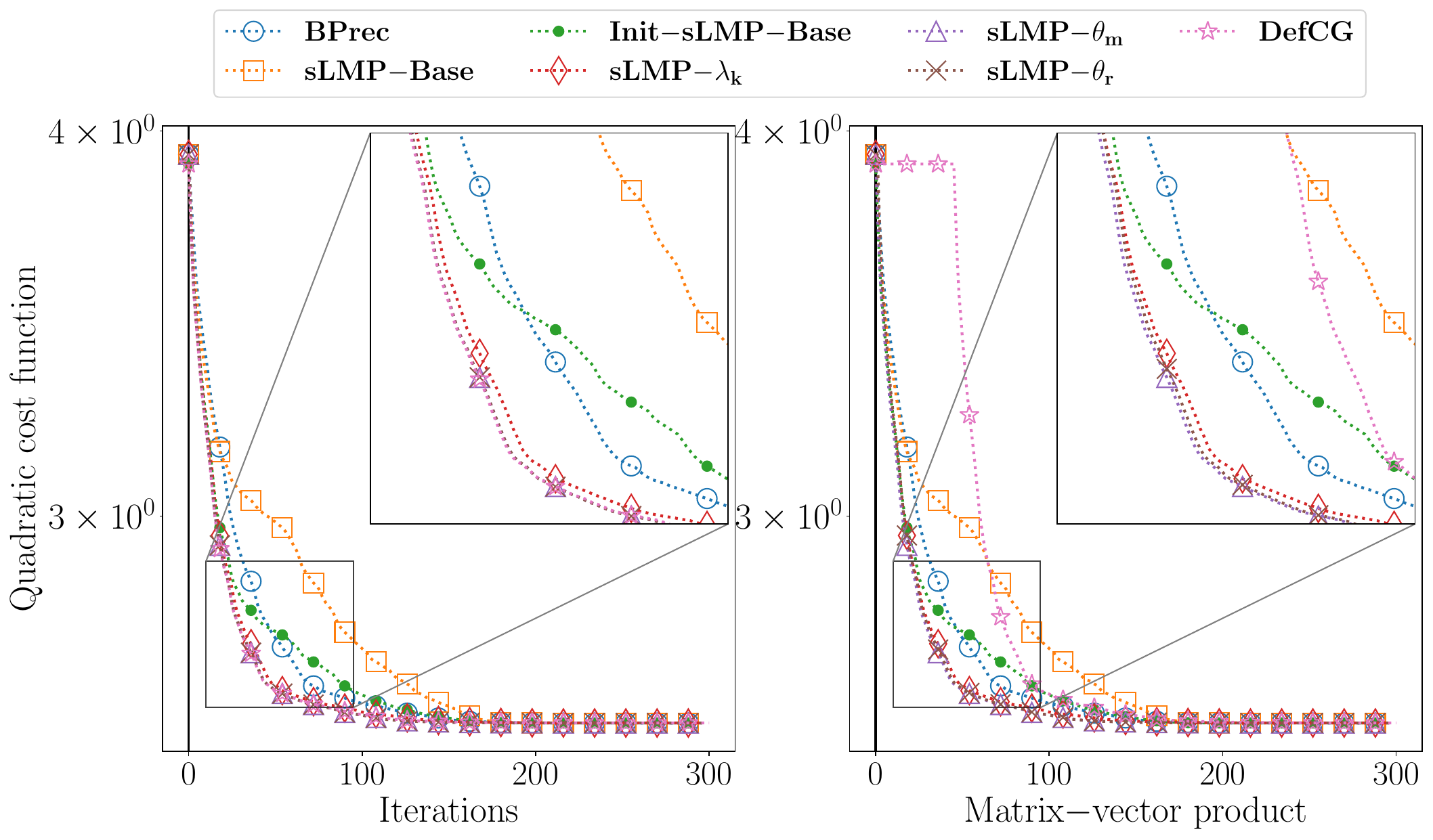}
    \caption{Quadratic cost function values along all CG iterates (left) and with respect to the number of matrix-vector product with the matrix $A^{(1)}$ and $A^{(2)}$ (right).}
    \label{fig:quadcost_CGiter}
\end{figure}

Note that, for \textbf{sLMP-$\theta_r$} we compute 
$\theta_r$ using~\eqref{eq:theta_r} 
with $r_0 = b^{(2)}$ and $A = A^{(1)}$. As a result, computation of approximate $\theta_r$ requires an extra matrix vector product with $A^{(1)}$. \Cref{fig:quadcost_CGiter} shows the quadratic cost function values~\eqref{Prob:quadratic} and number of matrix-vector products with $A^{(1)}$ and $A^{(2)}$ along CG iterations. 
\begin{figure}[ht]
    \centering
        \includegraphics[width=.9\linewidth]{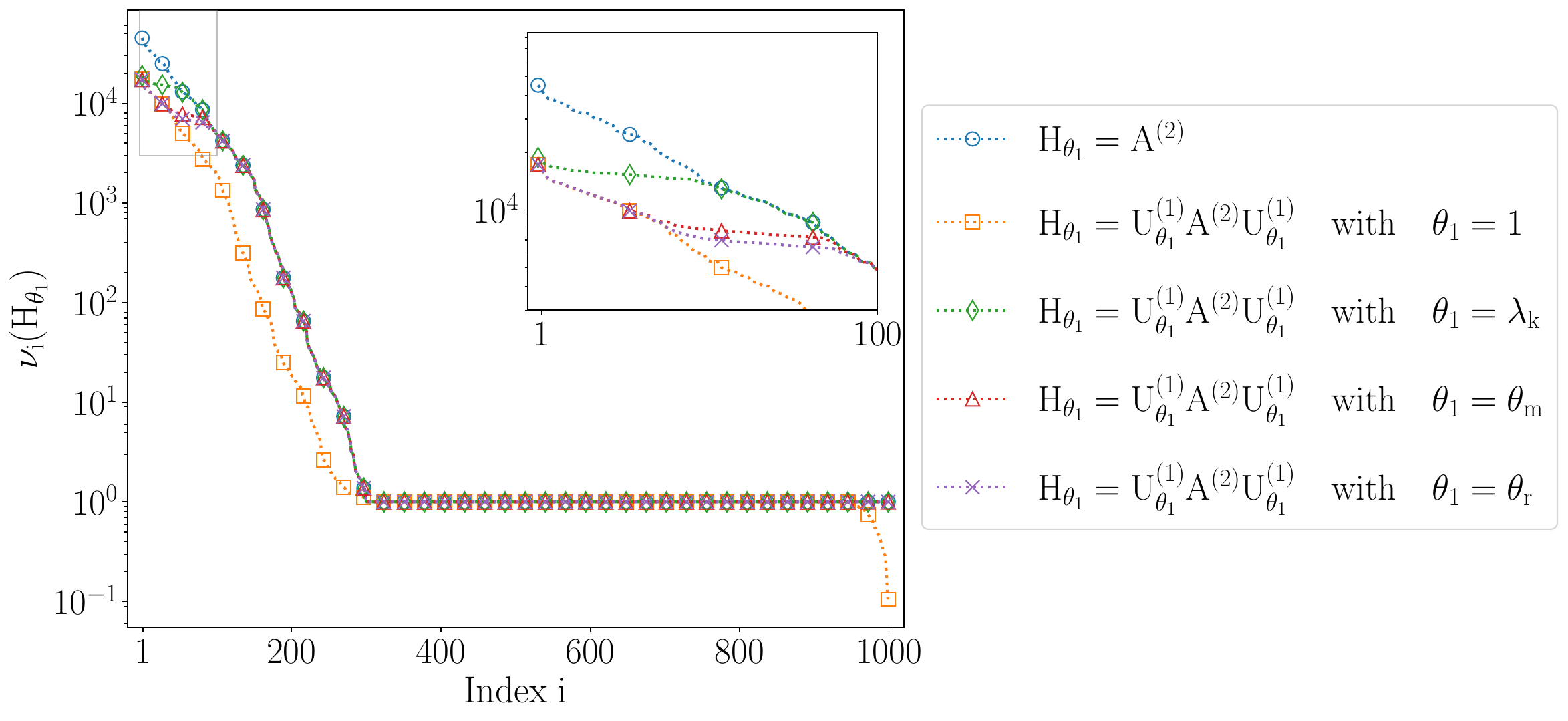}
    \caption{Spectrum of \(U_{\theta_1}^{(1)} A^{(2)} U_{\theta_1}^{(1)}\) for different values of \(\theta_1\) on a logarithmic scale. LowObs scenario $(k = 45)$.}
    \label{fig:spectrum}
\end{figure}

We can easily see that \textbf{sLMP-Base} is not necessarily better than \textbf{BPrec} especially in the early iterations. This means that the scaled spectral LMP, clustering the largest $k$ eigenvalues around $1$, might reduce the total number iterations to converge, however it does not guarantee better convergence for early iterations. The slow convergence of \textbf{sLMP-Base} can be partly explained by the fact that perturbations may cause some eigenvalues to appear near zero, as depicted in~\Cref{fig:spectrum}. When changing the clustering position from $1$
to $\lambda_k$ by using \textbf{sLMP-$\lambda_k$}, we can see that the method performs better than \textbf{BPrec}. In this case, however the gap between the cluster and the remaining spectrum as defined in~\Cref{theorem:midrange}, i.e. 
$
\alpha(\theta_1^{(1)}) / \theta_1^{(1)}
$,
can be large. When clustering around $\theta_r$ and $\theta_m$ is applied with \textbf{sLMP-$\theta_r$} and \textbf{sLMP-$\theta_m$} respectively, the value of $
\alpha(\theta_1^{(1)}) / \theta_1^{(1)}$ reduces for both cases (see Fig.~\ref{fig:spectrum}). This improves the convergence compared to \textbf{sLMP-$\lambda_k$} as seen from ~\Cref{fig:quadcost_CGiter}.

\textbf{Init-sLMP-Base} performs better than \textbf{sLMP-Base}, i.e. starting from \(x_0= S_k \Lambda_k^{-1}S_k^\top b^{(2)}\) improves performance compared to starting from $x_0 = 0$. This improvement arises because the initial residual's components in the eigenbasis of $A^{(2)}$ are reduced. In fact, without any perturbation, these components would be completely eliminated. Although, the performance is improved with this initial guess, it can not reach the performance of \textbf{DefCG}. This demonstrates that modifying the initial guess enhances convergence; however, the placement of the eigenvalue clustering can have an even more significant impact. This is evident from the fact that the performance of \textbf{sLMP-$\theta_m$} and \textbf{sLMP-$\theta_r$}  are very close to that of \textbf{DefCG}.

The right panel of~\Cref{fig:quadcost_CGiter} shows the values of the quadratic cost function as a function of the number of matrix-vector products performed with $A^{(j)}$ for $j=1,2$ across different methods. 
Although \textbf{DefCG} performs better, it is computationally expensive as it requires forming the projected matrix \(S_k^\top A^{(2)} S_k\).
Among the other techniques, \textbf{sLMP-$\theta_r$} requires one additional matrix-vector product with $A^{(1)}$ to compute $\theta_r$. However, as shown in~\Cref{fig:quadcost_CGiter}, \textbf{sLMP-$\theta_m$} and \textbf{sLMP-${\lambda_k}$} do not require any extra matrix-vector products either $A^{(1)}$ or $A^{(2)}$. 

These results indicate that the performance of CG, when used with scaled spectral LMP, can be significantly improved, approaching that of deflated CG, by selecting the position of the eigenvalue clusters based on CG's convergence properties. The cluster position is determined by $\theta$, whose computation incurs no additional cost for \textbf{sLMP-$\theta_m$} 
and \textbf{sLMP-${\lambda_k}$}.
Conclusions from experiments with \textit{HighObs} are very similar, the obtained results are depicted in \Cref{fig:quadcost_CGiterHighObs,fig:spectrumHighObs} in \Cref{appendix:GraphesHighObs}.

\section{Conclusion}
\label{sec:conclusion}

We have proposed a \textit{scaled} spectral LMP to accelerate the solution of a sequence of SPD systems $A^{(j)} x^{(j)} = b^{(j)}$ for $j \ge1 $. The \textit{scaled} LMP incorporates a low-rank update based on $k$ eigenpairs of the matrix $A$. We have provided theoretical analysis of the \textit{scaled} spectral LMP when $A^{(j)} = A$. We have shown that the scaled spectral LMP~\eqref{eq:ftheta} clusters $k$ eigenvalues around the scaling parameter $\theta$, and leaves the rest of the spectrum untouched. 

We have focused on the choice of $\theta$ to ensure that PCG achieves faster convergence, particularly in the early iterations. In the first approach, we have proposed choosing $\theta$ to guarantee a lower energy norm of the error at each iteration of PCG. In the second approach, we have obtained an optimum $\theta$ in the sense that it minimizes the energy norm of the error at the first iteration. Our analysis reveals that, with the optimal $\theta$, the components of the first residual is eliminated from the eigenspace of $A$, which aligns with the core principle of deflated CG. Lastly, we have also explored a scaling parameter that approximates the iterates of deflated CG. We have provided the link between the deflated CG and PCG with the scaled spectral LMP. 

We have compared different methods for solving a nonlinear weighted least-squares problem arising in data assimilation. In our numerical experiments, we used approximate eigenpairs to construct the scaled spectral LMP. First, we have demonstrated that selecting $\theta$ based on PCG convergence properties significantly accelerates early convergence compared to the conventional choice of $\theta=1$. Then, we have shown that $\theta$ values that reduce the spectral gap between $\theta$ and the remaining eigenvalues lead to faster convergence. 
Additionally, we have compared the scaled spectral LMP with deflated CG, showing that the scaled spectral LMP produces iterates similar to deflated CG, but at a negligible computational cost and memory, unlike deflated CG.
These numerical results clearly highlight the importance of selecting the preconditioner not only as an approximation to the inverse of $A$, but also with consideration of its role within PCG. In particular, we have demonstrated the significance of the placement of clustered eigenvalues, an often overlooked factor in the literature, on the early convergence of PCG.

As the next step, we will provide a detailed theoretical perturbation analysis in a forthcoming paper. Additionally, we aim to validate the proposed preconditioner in an operational weather prediction system. 

\appendix
\section{Deflated CG with \boldmath{\texorpdfstring{$S_k$}{Sk}}} 
\label{appendix:DefCG}
The deflation technique outlined in \Cref{algo:deflation} is defined for any deflation subspace $W$, see~\cite{DeflatedCG} for more details. The main idea is to speed-up the CG starting from an initial point such that the initial residual does not have components in the deflation subspace $W$ and to update the search directions such that $W^\top A p_j = 0$. A widely used approach is to choose $W$ as the eigenvectors corresponding to the eigenvalues that slows down the CG convergence.
\begin{algorithm}[ht]
\caption{Deflated-CG}\label{al:Deflated-CG}
\begin{algorithmic}[1]
\STATE Choose $k$ linearly independent vectors $w_1, w_2, \ldots, w_k$.
\STATE Define $W = [w_1, w_2, \ldots, w_k]$, and choose $x_{-1}$.
            \STATE Set $x_0^{\text{Def}} = x_{-1} + W(W^\top A W)^{-1} W^\top r_{-1}$, where $r_{-1} = b - A x_{-1}$. \COMMENT{$W^\top r_0 = 0$}
            \STATE Set $p_0  = r_0 - W \left(W^\top A W\right)^{-1} W^\top A r_0$. \COMMENT{$W^\top Ap_0 = 0$}
            \FOR{$j = 1, 2, \ldots$}
                \STATE $\alpha_{j-1} = r_{j-1}^\top r_{j-1} / (p_{j-1}^\top A p_{j-1})$
                \STATE $x_j^{\text{Def}} = x_{j-1}^{\text{Def}} + \alpha_{j-1} p_{j-1}$
                \STATE $r_j = r_{j-1} - \alpha_{j-1} A p_{j-1}$ \COMMENT{$W^\top r_j = 0$}
                \STATE $\beta_{j-1} = r_j^\top r_j / (r_{j-1}^\top r_{j-1})$
                \STATE $p_j = \beta_{j-1} p_{j-1} + r_j - W \left(W^\top A W\right)^{-1} W^\top A r_j$ \COMMENT{$W^\top Ap_j = 0$}
            \ENDFOR 
\end{algorithmic}
\label{algo:deflation}
\end{algorithm}

If we choose $W = S_k$, and using the fact that $S_k^\top A S_k = \Lambda_k$ and $ A S_k = S_k \Lambda_k$, we can achieve the following simplifications:
\begin{itemize}
    \item $x_0^{\text{Def}} = x_{-1} + S_k \Lambda_k^{-1} S_k^\top r_{-1}$, 
    \item $p_0 = r_0 - S_k S_k^\top r_0$.
    \item $p_j = \beta_{j-1} p_{j-1} + r_j - S_k S_k^\top r_j$.
\end{itemize}
\begin{lemma}
\label{lemma:deflatedCG}
The residual $r_j$ and the direction $p_j$ are orthogonal to $\text{span}(S_k)$.
\end{lemma}

\begin{proof}
We proceed by induction.
For $j = 0$,  $r_0 = r_{-1} - S_k S_k^{\top} r_{-1}$,
from which it follows that $S_k^{\top} r_0 = 0$. As a consequence, $S_k^{\top} p_0 = 0$.
Assume that $r_j$ and $p_j$ are orthogonal to $\text{span}(S_k)$ for $j$. 
We have  $r_{j+1} = r_{j} - \alpha_{j} A p_{j}$. From~\cite[Proposition 3.3]{DeflatedCG}, replacing $W$ by $S_k$, we have $S_k^T A p_{j} =0 $. Since $p_j , r_j \perp \text{span}(S_k)$ by assumption, it follows that $r_{j+1} \perp \text{span}(S_k)$.
For $p_{j+1} = \beta_{j} p_{j} + r_{j+1} - S_k S_k^\top r_{j+1} = \beta_{j} p_{j} + r_{j+1}$, we get $p_{j+1} \perp \text{Span}(S_k)$ since $S_k^\top r_{j+1} = 0$ as shown and $p_{j} \perp \text{Span}(S_k)$ by assumption.
\end{proof}
From \Cref{lemma:deflatedCG}, it follows that 
\(
p_j = \beta_{j-1} p_{j-1} + r_j - S_k S_k^\top r_j = \beta_{j-1} p_{j-1} + r_j.
\)
With these simplifications, it is clear that in exact arithmetic, deflated CG, when used with the deflated subspace consisting of a set of eigenvectors of A, generates iterates equivalent to those generated by using the initial guess $x_0^{Def}$ in standard CG.

\section{Results for the \textit{HighObs} scenario} 
\label{appendix:GraphesHighObs}
\small

\begin{figure}[ht]
\centering
    \includegraphics[width=.8\linewidth]{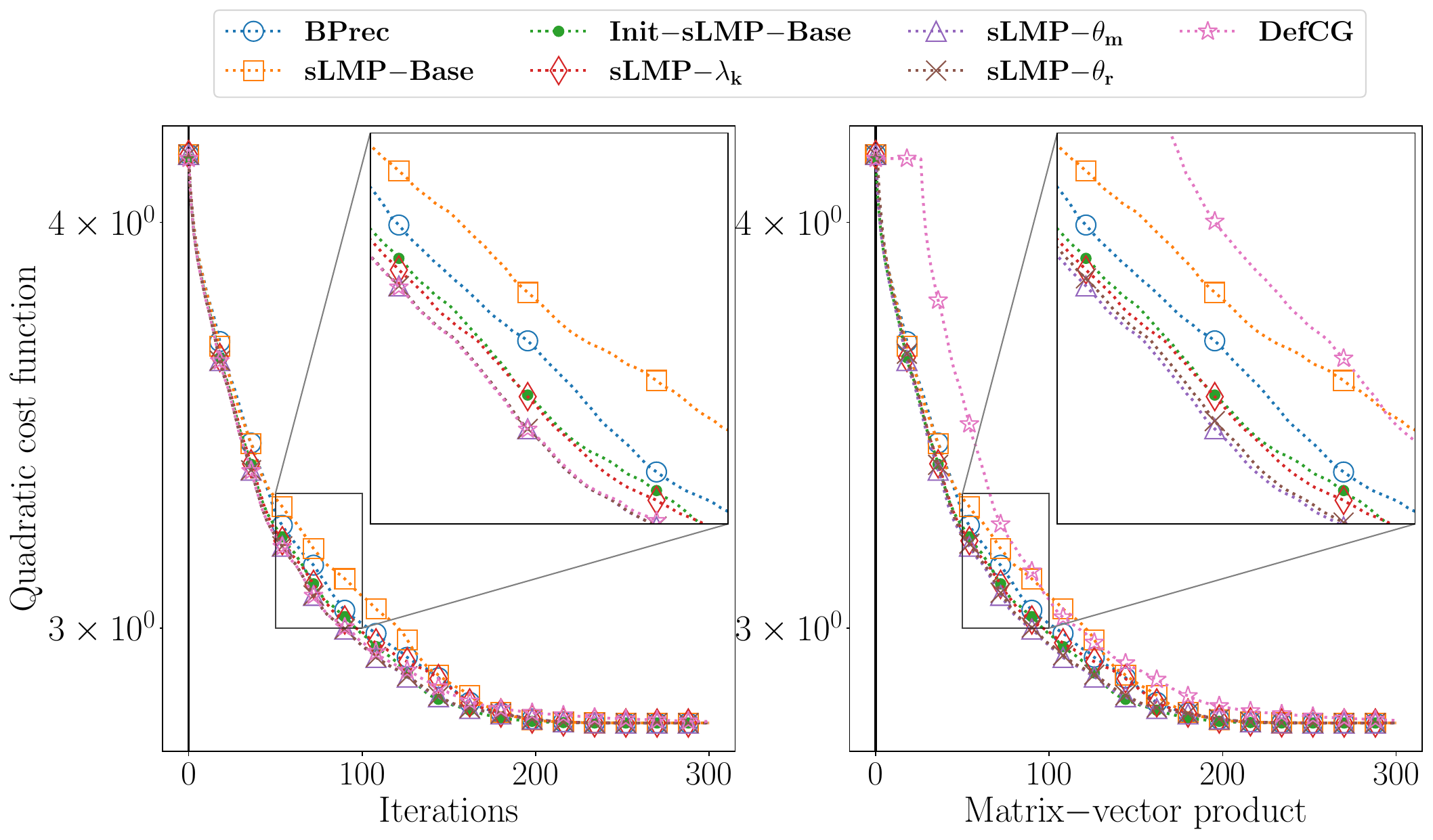}
    \caption{Quadratic cost function values along all CG iterates and with respect to the number of matrix-vector product for the \textit{HighObs} scenario.}
    \label{fig:quadcost_CGiterHighObs}
\end{figure}
\begin{figure}[ht]
    \centering
        \includegraphics[width=.8\linewidth]{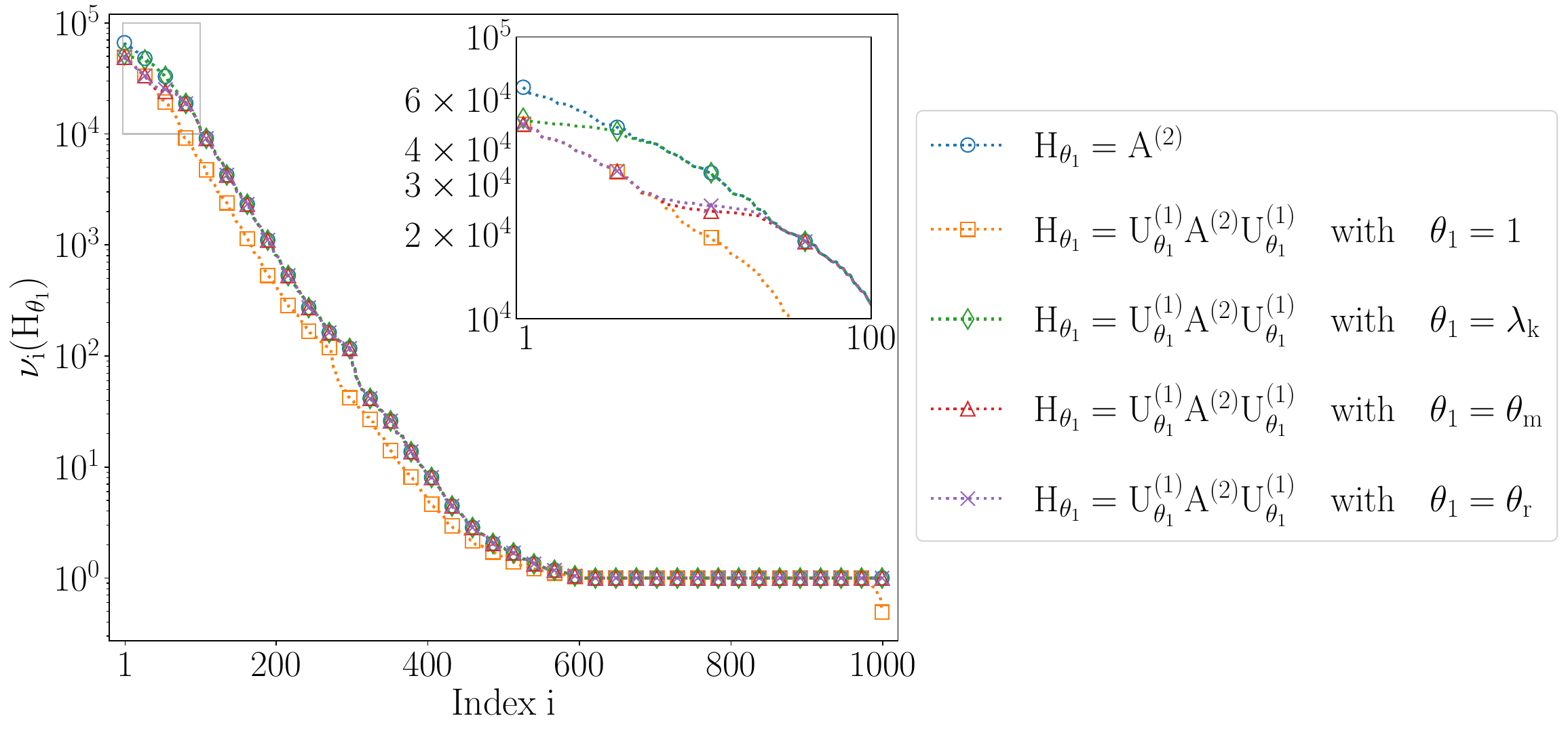}
    \caption{Spectrum of \(U_{\theta_1}^{(1)} A^{(2)} U_{\theta_1}^{(1)}\) with different \(\theta_1\) for the \textit{HighObs} scenario $(k = 26)$.}
    \label{fig:spectrumHighObs}
\end{figure}

\bibliographystyle{siamplain}
\bibliography{references}
\end{document}